\documentclass[11pt,letterpaper]{amsart}
\usepackage[T1]{fontenc}
\usepackage[utf8]{inputenc}
\usepackage{lmodern,amsmath,amssymb,amsthm,mathtools}
\usepackage{microtype,xcolor}
\usepackage[colorlinks=true,linkcolor=blue!55!black,
  citecolor=blue!55!black,urlcolor=blue!55!black]{hyperref}

\numberwithin{equation}{section}
\theoremstyle{plain}
\newtheorem{theorem}{Theorem}[section]
\newtheorem{proposition}[theorem]{Proposition}
\newtheorem{lemma}[theorem]{Lemma}
\newtheorem{corollary}[theorem]{Corollary}

\theoremstyle{definition}

\theoremstyle{remark}
\newtheorem{remark}[theorem]{Remark}

\newcommand{\Ric}{\operatorname{Ric}}
\newcommand{\Div}{\operatorname{div}_g}
\newcommand{\Vol}{\operatorname{Vol}}
\newcommand{\mr}{\mathbb{R}}
\newcommand{\ud}{\mathrm{d}}
\newcommand{\dmu}{\mathrm{d}\mu_g}
\newcommand{\dmuh}{\mathrm{d}\mu_h}
\newcommand{\Hess}{\operatorname{Hess}}
\newcommand{\ms}{\mathbb{S}}
\newcommand{\ip}[2]{\langle #1,#2\rangle}
\newcommand{\eps}{\varepsilon}

\allowdisplaybreaks[2]

\title[Fundamental group and  $Q$-curvature]
  {Fundamental group and $Q$-curvature on open four-manifolds}
\author{Mingxiang Li}
\address{Department of Applied Mathematics, The Hong Kong Polytechnic University,
  Hong Kong, China}
\email{mingxiang.li@polyu.edu.hk}
\date{}

\begin{document}

\begin{abstract}
Let $(M,g)$ be  a complete, connected, and non-compact four-manifold with nonnegative  scalar curvature and   nonnegative $Q$-curvature.  If $g$ is not Ricci-flat,  we show that  both the  fundamental group $\pi_1(M)$ and the  integral of $Q$-curvature are finite, satisfying
$$|\pi_1(M)|\int_MQ_g\ud\mu_g\leq 8\pi^2$$
with  the equality holding  if and only if $(M,g)$ is isometric to   $(\mathbb{R}^4, e^{2u}|\mathrm{d}x|^2)$ with vanishing isoperimetric ratio.
\end{abstract}

\maketitle

\section{Introduction}

Curvature conditions impose strong restrictions on the topology of
complete Riemannian manifolds.   In dimension two, this relationship is
expressed by the Gauss-Bonnet formula for a closed surface
\begin{equation*}
  \int_{M^2}K_g\,\dmu=2\pi\chi(M^2),
\end{equation*}
and by the Cohn-Vossen inequality for an open surface
with nonnegative Gaussian curvature
\begin{equation}\label{Cohn-Vossen}
  \int_{M^2}K_g\,\dmu\leq 2\pi\chi(M^2).
\end{equation}
Here $K_g$ is  the  Gaussian curvature and $\chi(M^2)$ is the Euler
characteristic. 

In higher dimensions,  scalar curvature provides one extension of this
picture. The famous Geroch conjecture asserts that a torus  admits no complete  metric of positive scalar
curvature, proved by Schoen-Yau in dimensions three through seven and
by Gromov-Lawson in all dimensions. Complete metrics with positive scalar curvature on noncompact manifolds also satisfy substantial topological restrictions; see \cite{LUY}.

  $Q$-curvature is another natural generalization  of the Gaussian curvature in   dimension four.   On the one hand, Branson's $Q$-curvature is given by
\begin{equation}\label{def-of-Q}
  Q_g=\frac16\bigl(-\Delta_gR_g-3|\Ric_g|_g^2+R_g^2\bigr)
\end{equation}
where $R_g, \Ric_g, \Delta_g$ denote the scalar curvature, the Ricci tensor, and  the Laplace-Beltrami operator respectively.
On the other hand, the  Paneitz operator is given by
\begin{equation}\label{Paneitz-operator}
  P_gu=\Delta_g^2u-\Div\left[
    \left(\frac23R_gg-2\Ric_g\right)^{\sharp}(\nabla u)\right]
\end{equation}
where the notation $T^{\sharp}$ means that
one index is raised using $g$.
They are linked by the beautiful conformal covariant equation
\begin{equation}\label{conformal-Q}
  P_gu+Q_g=Q_{e^{2u}g}e^{4u}.
\end{equation}
In particular, the integral of $Q$-curvature is conformally invariant on closed four-manifolds. We refer interested readers to Chang’s excellent lecture note \cite{Chang} for an overview of $Q$-curvature on four-manifolds.

 For a closed four-manifold $(M^4,g)$, the Chern-Gauss-Bonnet formula
takes the form
\begin{equation*}
  \int_{M^4}\left(\frac14|W_g|_g^2+Q_g\right)\dmu
  =8\pi^2\chi(M^4),
\end{equation*}
where $W_g$ denotes the  Weyl tensor. It plays a fundamental role in studying geometric and topological properties of four-manifolds.
For closed, oriented four-manifolds, 
the  Hirzebruch  signature formula is given by
\[
  \tau(M^4)=\frac{1}{48\pi^2}
  \int_{M^4}\bigl(|W_g^+|_g^2-|W_g^-|_g^2\bigr)\dmu
\]
where $W_g^{\pm}$denote the  self-dual and anti-self-dual component of $W_g$. Combining these two formulas shows that positive $Q$-curvature everywhere immediately implies the following topological restriction
\begin{equation*}
	Q_g>0\Longrightarrow 2\chi(M^4)>3|\tau(M^4)|.
\end{equation*}
See Gursky’s remarkable work  \cite{Gursky Ann} for  more related results and deep discussions.  On the torus, Lin and Yuan \cite{LinYuan} studied the Geroch-type conjecture for $Q$-curvature. 

Throughout this paper,   all the Riemannian  manifolds are assumed to be smooth, connected, and without boundary. An \emph{open  (closed) manifold}  is a complete and noncompact (compact)  Riemannian  manifold.  

Our first result shows that positive
$Q$-curvature alone imposes no topological obstruction on open four-manifolds which differs from the closed case. 

\begin{theorem}\label{thm:main1}
Let $(M,g_0)$ be an open four-manifold. There exists a complete  metric $g$ on $M$ such that $Q_g>0$ everywhere.
\end{theorem}

Therefore, additional curvature conditions need to be taken into account, and the scalar curvature is a very natural candidate. Let $(M,g)$ be  a closed  four-manifold. If  the $Q$-curvature and the scalar curvature satisfy
\begin{equation}\label{QRpositive}
	Q_g>0 \quad \mathrm{and} \quad R_g>0,
\end{equation}
Chang-Gursky-Yang \cite{CGY} showed that there exists  a conformal metric $\tilde g\in[g]$ such that \begin{equation}\label{positive-Ricci}
	0<\Ric_{\tilde g}<\frac12R_{\tilde g}\tilde g.
\end{equation}
Then the compactness of $M$, together with the Bonnet-Myers theorem, implies  the fundamental group $\pi_1(M)$  is finite.  Assuming $R_g\geq 0$,
Gursky \cite{Gursky CMP}  showed that 
\begin{equation}\label{Gursky-inequlaity}
	\int_M Q_g\dmu\leq 16\pi^2
\end{equation}
with   the equality holding if and only if  $(M,g)$ is conformally equivalent to $(\ms^4,g_{\ms^4})$.
Combining \eqref{positive-Ricci} with \eqref{Gursky-inequlaity} (see \cite[Proposition 5.1]{JLW}), under  \eqref{QRpositive} on the  closed manifold, one has 
\begin{equation}\label{CGY-inequality}
	|\pi_1(M)|\int_MQ_g\ud\mu_g\leq 16\pi^2.
\end{equation}

Recently, a Bonnet-Myers type theorem for $Q$-curvature on complete four-manifolds has been established  by Jiang-Li-Wang \cite{JLW}, where  we showed that if a complete four-manifold $(M,g)$ satisfies 
$$Q_g\geq 1 \quad \mathrm{and}\quad R_g\geq 0,$$ then $M$ must be compact. Consequently, $\pi_1(M)$ is finite and \eqref{CGY-inequality} still holds.

Motivated by the work of Chang-Gursky-Yang \cite{CGY} and our previous work \cite{JLW}, it is natural to ask whether simultaneous nonnegativity of the $Q$-curvature and scalar curvature can effectively constrain the topology of open four-manifolds? If $g$ is Ricci-flat, then both $Q_g$ and $R_g$ vanish, and the control they provide reduces to the study of Ricci-flat manifolds which is a very grand  topic that has been extensively studied and is quite delicate. Therefore, we do not treat this case in this paper.  

Before stating our main theorem, we  introduce another geometric quantity.
For an open $n$-manifold $(M^n,g)$, the isoperimetric   ratio  $I_g(M^n)$ is defined as
$I_g(M^n)=\inf_{\Omega\subset M}\frac{\Vol_g(\partial\Omega)^{\frac{n}{n-1}}}{\Vol_g(\Omega)}$
where the infimum is taken over all non-empty smooth bounded domains
$\Omega\subset M^n$.

\begin{theorem}\label{main theorem}
	Let $(M^4,g)$ be an open four-manifold with $Q_g \geq 0$ and  $R_g \geq 0$.  Suppose that $g$ is not Ricci-flat.
	Then  both the fundamental group $\pi_1(M^4)$ and the integral of $Q$-curvature  are finite.  Moreover, 
	\begin{equation}\label{CV-type-inquality}
		|\pi_1(M^4)|\int_{M^4}Q_g\dmu\leq 8\pi^2
	\end{equation}
	with the equality holding if  and only if  $(M^4,g)$ is  isometric to $(\mathbb{R}^4, e^{2u}|\ud x|^2)$  for some smooth function $u$ with  isoperimetric    ratio $I_{e^{2u}|\ud x|^2}(\mr^4)=0$.
\end{theorem}

To better interpret   our topological result concerning  $Q$-curvature, we recall the corresponding restrictions imposed by sectional and Ricci curvatures in contrast.

The soul theorem of Cheeger-Gromoll \cite{Cheeger-Gromoll} implies that an open manifold with nonnegative sectional curvature is diffeomorphic to the normal bundle of a compact soul, and in particular has finitely generated fundamental group.  Gromoll-Meyer \cite{Gromoll-Meyer} showed that positive sectional curvature forces the open manifold to be diffeomorphic to $\mr^n$  and   Perelman \cite{Perelamn} proved the same conclusion when sectional curvature is nonnegative and is positive at one point.

When $\Ric_g\geq kg$ for some constant $k>0$, the famous Bonnet-Myers theorem shows that the manifold must be compact and then the fundamental group is finite. Sixty years ago, Milnor studied the relationship between curvature and the fundamental group and conjectured that, for $n\geq 3$, every complete 
$n$-dimensional Riemannian manifold with nonnegative Ricci curvature has a finitely generated fundamental group.
   For an open three-manifold, Schoen-Yau \cite{SY-82} proved that positive Ricci curvature forces the manifold to be diffeomorphic to $\mr^3$. Building on their work,  Liu \cite{Liu} established the three-dimensional case of Milnor's conjecture.   In higher dimensions, Bru\`e-Naber-Semola \cite{BNS} constructed counterexamples to Milnor's conjecture in dimensions at least seven, and their subsequent work \cite{BNS-JEMS} produced a six-dimensional counterexample. The recent work  of Wu-Yan \cite{WY} gives examples with positive Ricci curvature and infinitely generated fundamental group in dimensions four and five.

 The inequality \eqref{CV-type-inquality} can be viewed as a Cohn-Vossen-type inequality analogous to \eqref{Cohn-Vossen}. Let $(\mr^4,g=e^{2u}|\ud x|^2)$ be  a complete and conformally flat manifold. If  $R_g\geq0$ and $Q_g\in L^1(\mr^4,g)$, Chang, Qing, and Yang have proved  \eqref{CV-type-inquality} in \cite{CQY}. In \cite{LWX}, we  replaced the assumption $Q_g\in L^1(\mr^4,g)$ by the negative part $Q_g^-\in L^1(\mr^4,g)$.    Moreover,  if $Q_g\geq 0$ and $R_g\geq 0$,  we  \cite{LWX} showed that $\Ric_g\geq 0$.  Recently, we  \cite{LX} showed  that even the sectional curvature is nonnegative. These results reveal an interesting connection between the  $Q$-curvature and the Ricci and sectional curvatures under restrictive conformal structural assumptions.

Now we briefly outline the organization  of this paper. First, in Section \ref{sec:positive-Q}, we prove Theorem \ref{thm:main1} by using Morse theory and the conformal structure of the Paneitz operator.  Second, in Section \ref{sec:splitting}, assuming  $Q_g\geq 0$ and $R_g\geq 0$, we establish a splitting theorem for open four-manifolds.  Later,   in Section \ref{sec:green-function}, we recall some basic facts about the Green function of the conformal Laplacian operator  on open manifolds and establish some integral estimates involving such Green functions. Then,  in Section \ref{sec:HYinequality}, we show that  the integral of $Q$-curvature is at most $16\pi^2$; see Theorem \ref{thm:CVandHY-for-Q}. With the help of such  inequality and the splitting theorem,  in Section \ref{sec:pi1}, we prove  the fundamental group $\pi_1(M)$  is finite under $Q_g\geq 0$ and $R_g>0$. In Section \ref{sec:modified-Green}, we establish some  integral estimates for a potential of conformal Laplacian operator. Then, using these estimates, we improve the upper bound of  the integral of $Q$-curvature to $8\pi^2$ in Section \ref{sec:impoved-HY}. Finally, we complete the proof of Theorem \ref{main theorem} and discuss  the impact of $Q$-curvature in other dimensions in Section \ref{sec:proof}.

\section{Open four-manifolds with positive $Q$-curvature}
\label{sec:positive-Q}

Throughout this section, let $(M^4,g_0)$ be an open four-manifold. We shall construct a complete
metric $g\geq g_0$ with $Q_g>0$. The construction first modifies $g_0$
near the critical points of a proper Morse function and then makes a
conformal change of the modified metric.

The standard Morse theory for open manifolds provides a smooth proper Morse function
\[
f:M^4\longrightarrow [0,\infty)
\]
with no local maxima; see \cite[Lemma 1.1]{Phillips}. Its critical set is closed and
discrete. At each critical point, the Hessian of $f$ is nondegenerate and
has at least one positive eigenvalue. We first adjust the metric near these critical points.

\begin{lemma}\label{lem:metric-near-critical}
There exists a complete  metric $g_1\geq g_0$ such that
\begin{equation}\label{eq:adapted}
\Delta_{g_1}f+2\lambda_{\min}(\Hess_{g_1}f)>0
\end{equation}
at every critical point of $f$, where $\lambda_{\min}(\Hess_{g_1}f)$ denotes  the smallest eigenvalue  of $\Hess_{g_1}f$ with respect to $g_1$.
\end{lemma}

\begin{proof}
First, fix a critical point $p$. Since $\ud f(p)=0$, the Hessian at $p$
is independent of the connection. Diagonalize this bilinear form in an
orthonormal basis with respect to $g_0$, and denote its eigenvalues by
$\lambda_1,\ldots,\lambda_4$. Multiply $g_0$ by a constant $\Lambda\geq1$ on the negative eigenspaces. If all the eigenvalues are positive, \eqref{eq:adapted} already holds.
Otherwise, as $\Lambda\to\infty$, the left-hand side of
\eqref{eq:adapted} converges to the sum of the positive eigenvalues,
which is strictly positive.  Then, choose pairwise disjoint, locally finite, relatively compact coordinate neighborhoods of the critical points. In each neighborhood, extend the
prescribed positive semidefinite form smoothly as a positive
semidefinite tensor, and multiply it by a nonnegative smooth cutoff
equal to one near the critical point. Adding these locally finite
modifications to $g_0$ produces a complete metric $g_1\geq g_0$.  
\end{proof}

\begin{proof}[Proof of Theorem~\ref{thm:main1}]
Fix the metric $g_1$ obtained  from Lemma~\ref{lem:metric-near-critical}. Unless
otherwise specified, all derivatives and norms in this proof are taken
with respect to $g_1$. For any smooth function
$F:[0,\infty)\to\mr$, the chain rule and Bochner's formula give
\begin{equation}\label{eq:chain}
P_{g_1}(F(f))
=|\nabla f|^4F^{(4)}+\varphi_1F'''+\varphi_2F''+(P_{g_1}f)F',
\end{equation}
where the derivatives of $F$ are evaluated at $f$ and
\begin{equation}\label{eq:coefficients}
\begin{aligned}
\varphi_1
&=2|\nabla f|^2\Delta f+4\nabla^2f(\nabla f,\nabla f),\\
\varphi_2
&=(\Delta f)^2+2|\nabla^2f|^2
  +4\langle\nabla f,\nabla\Delta f\rangle\\
&\qquad+4\Ric_{g_1}(\nabla f,\nabla f)
  -\frac23R_{g_1}|\nabla f|^2.
\end{aligned}
\end{equation}

At a critical point $p$, one has 
\[
\varphi_1(p)=0,
\qquad
\varphi_2(p)=(\Delta f(p))^2+2|\nabla^2f(p)|^2>0.
\]
Moreover,
\[
\varphi_1\geq
2|\nabla f|^2\bigl(\Delta f+2\lambda_{\min}(\nabla^2f)\bigr).
\]
By \eqref{eq:adapted} and continuity, each critical point therefore has
a relatively compact neighborhood $U_p$ such that
$\varphi_1\geq0$ and $\varphi_2>0$ on the closure $\overline{U_p}$. These
neighborhoods may be chosen locally finite. Set
$U=\bigcup_pU_p$. Then
\begin{equation}\label{eq:regions}
\varphi_1\geq0,\quad \varphi_2>0\quad\hbox{on }U,
\qquad
|\nabla f|>0\quad\hbox{on }M\setminus U.
\end{equation}

For every integer $j\geq0$, the shell $f^{-1}([j,j+1])$ is compact.
Only finitely many $U_p$ meet this shell, and $\varphi_2$ is positive
on their compact closures. Hence $\varphi_2$ has a positive lower bound
on the intersection of the shell with $U$. The function $|\nabla f|$
also has a positive lower bound on  $f^{-1}([j,j+1]) \cap M\setminus U$.
Choose $\Lambda_j\geq1$ large enough that,
\begin{equation}\label{eq:bounds}
\begin{aligned}
\varphi_2\Lambda_j
&\geq |P_{g_1}f|+|Q_{g_1}|+1,
\quad \quad \quad \mathrm{on }\;  f^{-1}([j,j+1])\cap U,\\
|\nabla f|^4\Lambda_j
&\geq |\varphi_1|+|\varphi_2|+|P_{g_1}f|+|Q_{g_1}|+1,
\quad \mathrm{on }\;\;f^{-1}([j,j+1])\cap M\setminus U.
\end{aligned}
\end{equation}
Define a function
\[
\psi(t)=
\begin{cases}
\displaystyle
\frac{\int_0^t(t-s)e^{-1/s}\,\ud s}
     {\int_0^1(1-s)e^{-1/s}\,\ud s},&t>0,\\
0,&t\leq0.
\end{cases}
\]
It is not hard to check that the function $\psi$ is smooth, $\psi',\psi''\geq0$, and $\psi(1)=1$.
Set
\[
\kappa(t)=\Lambda_0+
\sum_{j=1}^{\infty}\Lambda_j\psi(t-j+1).
\]
Since the sum is locally finite, then $\kappa$ is smooth and satisfies
$\kappa\geq1$, $\kappa'\geq0$, and $\kappa''\geq0$. In addition, one has 
\begin{equation}\label{eq:kappa}
\kappa(t)\geq\Lambda_j\qquad(j\leq t\leq j+1).
\end{equation}

Choose 
\[
F(t)=\int_0^t
\exp\left(\int_0^s\kappa(r)\,\ud r\right)\ud s,
\qquad u=F(f).
\]
For $t\geq0$, we have $F(t)\geq t$, $F'\geq1$, and
\[
F''=\kappa F',\qquad
F'''=(\kappa^2+\kappa')F',\qquad
F^{(4)}=(\kappa^3+3\kappa\kappa'+\kappa'')F'.
\]
Consequently, one has
\begin{equation}\label{eq:derivatives}
F^{(4)}\geq\kappa F''',
\qquad F'''\geq F''\geq F'\geq1.
\end{equation}
In $U$, equations
\eqref{eq:chain}-\eqref{eq:kappa} give
\[
P_{g_1}u
\geq\varphi_2F''-|P_{g_1}f|F'
=(\varphi_2\kappa-|P_{g_1}f|)F'
\geq |Q_{g_1}|+1.
\]
In  $M\setminus U$, equations \eqref{eq:bounds} and
\eqref{eq:derivatives} yield
\begin{align*}
P_{g_1}u
&\geq |\nabla f|^4F^{(4)}
       -|\varphi_1|F'''-|\varphi_2|F''-|P_{g_1}f|F'\\
&\geq\bigl(|\nabla f|^4\kappa
       -|\varphi_1|-|\varphi_2|-|P_{g_1}f|\bigr)F'''\geq |Q_{g_1}|+1.
\end{align*}
Thus, there holds $$P_{g_1}u+Q_{g_1}\geq1, \quad \quad  \mathrm{in}\quad  M.$$
Choose  $g=e^{2u}g_1$.  The conformal property
\eqref{conformal-Q} yields 
\[
Q_g=e^{-4u}(P_{g_1}u+Q_{g_1})\geq e^{-4u}>0.
\]
Since $u\geq0$, we have $g\geq g_1\geq g_0$,
and hence $g$ is also complete.
\end{proof}

\section{A splitting theorem }
\label{sec:splitting}

The  following lemma is due to an observation of Gursky-Malchiodi \cite{GM}.
\begin{lemma}\label{lem:GM-lemma}
Let  $(M^4,g)$ be a complete four-manifold with  $Q_g\geq 0$ and $R_g\geq 0$. Then, either $\Ric_g\equiv 0$ or $R_g>0$.
\end{lemma}
\begin{proof}
	Using \eqref{def-of-Q}, $Q_g\geq 0$ and $R_g\geq 0$ yields
	$$-\Delta_g R_g	+R_g^2\geq 0.$$
	The strong maximum principle yields $R_g\equiv 0$ or $R_g>0$. When $R_g\equiv0$, inserting it  to \eqref{def-of-Q}, one has
	$$0\leq Q_g	=-\frac{1}{2}|\Ric_g|_g^2$$
	which yields $\Ric_g\equiv 0$.
\end{proof}
The algebraic structure of $Q$-curvature yields a lower bound for Ricci
curvature in the spectral sense observed by Chang-Gursky-Yang \cite{CGY}.

\begin{lemma}\label{lem:lambda_1}
Let  $(M^4,g)$ be a complete four-manifold. Set
\begin{equation}\label{def-of-rho}
	\rho(x)=\min_{\substack{v\in T_xM\\g(v,v)=1}}\Ric_g(v,v).
\end{equation} 
Then, there holds
\begin{equation*}
\left(-\frac12\Delta_g+\rho\right)R_g
\geq 3Q_g+2\rho^2.
\end{equation*}
\end{lemma}

\begin{proof}
At a fixed point, choose a $g$-orthonormal basis for the Ricci
tensor and denote  its eigenvalues as
$\rho=\lambda_1\leq\lambda_2\leq\lambda_3\leq\lambda_4$.
Cauchy's inequality yields
\begin{align*}
|\Ric_g|_g^2
&=\rho^2+\lambda_2^2+\lambda_3^2+\lambda_4^2\\
&\geq\rho^2+\frac{(R_g-\rho)^2}{3}
=\frac{R_g^2-2\rho R_g+4\rho^2}{3}.
\end{align*}
Inserting   \eqref{def-of-Q}, we obtain
\[
-\frac12\Delta_gR_g+\rho R_g
=3Q_g+\frac32|\Ric_g|_g^2-\frac12R_g^2+\rho R_g
\geq3Q_g+2\rho^2.
\]
\end{proof}

The preceding lemma implies that the operator $-\frac{1}{2}\Delta_g+\rho$ is nonnegative when $Q_g\geq 0$ and $R_g>0$. Combining with the spectral splitting results of Antonelli-Pozzetta-Xu \cite{APX} and Catino-Mari-Mastrolia-Roncoroni \cite{CMMR}, we  obtain the following splitting theorem for $Q$-curvature.

\begin{theorem}\label{thm:splitting}
	Let $(M^4,g)$ be an open four-manifold with $Q_g\geq0$ and $R_g\geq0$.
	Suppose that  $M^4$ has more than one end. Then it splits isometrically as
	\[
	(M^4,g)\cong(\mr\times N^3,\ud t^2\oplus g_N),
	\]
	where $(N^3,g_N)$ is a compact three-dimensional Einstein manifold with
	nonnegative Ricci curvature. Moreover, $Q_g\equiv0$.
\end{theorem}
\begin{proof}
First,  we claim
\begin{equation}\label{eq:spectral-nonnegative}
\int_M\left(\frac12|\nabla\zeta|^2+\rho\zeta^2\right)\dmu\geq0,
\qquad\forall \zeta\in C_c^\infty(M)
\end{equation}
where $\rho$ comes from \eqref{def-of-rho}.

 Lemma \ref{lem:GM-lemma} yields that either $
\Ric_g\equiv 0$ or $R_g>0$. If  $\Ric_g\equiv0$, then
\eqref{eq:spectral-nonnegative} follows immediately.

Suppose now that $R_g>0$. Integrating  by parts and then using 
Lemma \ref{lem:lambda_1}, we have
\begin{align*}
&\int_M\left(\frac12|\nabla\zeta|^2+\rho\zeta^2\right)\dmu\\
=&\frac12\int_M R_g^2
  \left|\nabla\left(\frac{\zeta}{R_g}\right)\right|^2\dmu+\int_M
\frac{(-\frac12\Delta_g+\rho)R_g}{R_g}\zeta^2\dmu\geq 0
\end{align*}
which proves \eqref{eq:spectral-nonnegative}.

Then, applying  the spectral splitting theorem of Antonelli-Pozzetta-Xu
\cite[Theorem 1.1]{APX}  or 
Catino-Mari-Mastrolia-Roncoroni \cite[Theorem 1.1]{CMMR}, we obtain
\[
(M^4,g)\cong(\mr\times N^3,\ud t^2\oplus g_N),
\]
where $(N^3,g_N)$ is compact and $\Ric_{g_N}\geq0$. For this product, $R_g=R_{g_N}$ and 
\begin{equation}\label{eq:product-q}
6Q_g=-\Delta_{g_N}R_{g_N}
-3\left|\Ric_{g_N}-\frac{R_{g_N}}3g_N\right|_{g_N}^2.
\end{equation}
Integrating  \eqref{eq:product-q} over the compact  manifold $N^3$ yields
\[
0\leq6\int_N Q_g\,\ud\mu_{g_N}
=-3\int_N\left|\Ric_{g_N}-\frac{R_{g_N}}3g_N\right|_{g_N}^2
\ud\mu_{g_N}\leq0.
\]
Consequently, one has $Q_g\equiv0$ and
$\Ric_{g_N}\equiv \frac{R_{g_N}}3g_N\geq 0$. Thus, we finish the proof.
\end{proof}

Using the above splitting theorem, we are able to  classify the locally conformally flat open manifolds with $Q_g\geq 0$ and $R_g>0$. 

\begin{proposition}\label{prop:cylinder}
	Let $(X^4,g)$  be a  simply connected and locally conformally flat open four-manifold with $Q_g\geq 0$ and $R_g>0$.  
	\begin{enumerate}
		\item 	If $Q_g\not \equiv 0$,  	then 
		$$(X^4,g)\cong (\mr^4, e^{2u}|dx|^2)$$ for some smooth function $u$.
		\item	If $Q_g\equiv 0$, then
		\[
		(X^4,g)\cong
		\bigl(\mr\times\ms^3,\ud t^2\oplus a^2g_{\ms^3}\bigr)
		\]
		where $a$ is a positive constant and $g_{\ms^3}$ is the round metric.
	\end{enumerate}
\end{proposition}

\begin{proof}
	Using Schoen-Yau's work  \cite{SY88} and Chodosh-Li's result \cite[Corollary 4]{CL} , we identify \(X\) conformally with
	$\Omega=\Phi(X)\subset\ms^4$ by the developing map $\Phi$. Then $\Phi^*(g)=e^{2w}g_{\ms^4}$ on $\Omega$.
	Theorem 1.5 of Ma-Qing \cite{MQ}  implies that $\ms^4\setminus\Omega$ consists of only finitely many
	points. Since $X$ is noncompact, $\ms^4\setminus \Omega$ is nonempty.
	
	If  $\ms^4\setminus\Omega$ consists of one point, the stereographic projection
	identifies $g$ with $e^{2u}|\ud x|^2$ on $\mr^4$ for some smooth function $u$.
	
	If $Q_g\not \equiv 0$, Theorem  \ref{thm:splitting} shows that  $\ms^4\setminus\Omega$ consists of one point. Thus, $(X^4,g)$ is isometric to  $(\mr^4, e^{2u}|dx|^2)$.
	
	If $Q_g\equiv 0$, $R_g\geq 0$ and $\ms^4\setminus \Omega$ consists of one point, then  one has $\Delta^2 u=0$ and $\Delta u+|\nabla u|^2\leq 0$. In particular, $\Delta u \leq 0$ and  then Liouville theorem yields $-\Delta u$ is equal  a nonnegative constant $c$ and then $|\nabla u|^2\leq -\Delta u=c$. Then
	$$c|B_r(0)|=\int_{B_r(0)}(-\Delta u)\ud x=-\int_{\partial B_r(0)}\partial_ru\ud\sigma\leq |\partial B_r(0)|\sqrt{c}.$$
	Letting $r\to\infty$, we show that $c=0$. Therefore, $u$ must be a constant and then $R_g\equiv 0$ contradicting the assumption $R_g>0$.
	Thus,  if $Q_g\equiv 0$, $\ms^4\setminus  \Omega$ consists of at least two points, so $X^4$ has at
	least two ends.  Using Theorem \ref{thm:splitting} and the simply connected assumption, we finish  the proof.
\end{proof}

\section{ Green function of conformal Laplacian operator}\label{sec:green-function}

Let $(M^n,g)$ be a complete manifold with dimension $n\geq 3$ throughout this section.
For  the closed case, it is well-known that   (see \cite{Aubin}) if $R_g>0$,  there exists  a positive Green function of  the conformal Laplacian operator.
 For the  open case,    the positive minimal Green function for some nonnegative elliptic operators has also been well-studied.  In particular,   if  $V$ is a nonnegative and non-zero smooth function,  the Schr\"odinger operator 
\begin{equation}\label{operator-L}
	Lu=-\Delta_gu +Vu
\end{equation}
 has a positive minimal Green function.  

For the reader's convenience, we sketch  the proof of the existence of a  positive, minimal, and  symmetric  Green function for such operator and establish some pointwise estimates for such Green function.  More details and discussions can be found in  \cite{Pinchover}.

\begin{lemma}\label{lem:green}
Let $(M^n,g)$ be an open manifold of
dimension $n\ge3$. Suppose that $V\in C^\infty(M)$ satisfies $V\ge0$ and
$V\not\equiv0$. For each  point $p\in M$, the operator $L$ defined in 
 \eqref{operator-L} has a  minimal positive symmetric Green function
$G_L(p,\cdot)$ satisfying $$LG_L(p,\cdot)=\delta_p(\cdot).$$
Moreover, for every smooth relatively compact neighborhood $U$ of $p$, one has 
\begin{equation}\label{eq:exterior}
 \sup_{x\in M\setminus U}G_L(p,x)
 \le \max_{x\in\partial U}G_L(p,x)<\infty.
\end{equation}
\end{lemma}

\begin{proof}
First, choose a smooth connected exhaustion
\[
 p\in\Omega_1\Subset\Omega_2\Subset\cdots,
 \qquad \bigcup_{j=1}^\infty\Omega_j=M.
\]
Since $V\ge0$, from the  standard elliptic theory (see \cite{Aubin}), there exists  a  positive symmetric Dirichlet Green function  $G_{L,\Omega_j}$  of the operator $L$ on $\Omega_j$ satisfying
$$LG_{L,\Omega_j}=\delta_p, \quad  \mathrm{in}\; \Omega_j;  \quad G_{L,\Omega_j}=0, \quad \mathrm{on} \; \partial \Omega_j.$$ 
For brevity, set $G_j=G_{L,\Omega_j}$. The maximum
principle gives
\[
 0<G_j\le G_{j+1}\quad\text{on }\Omega_j\setminus\{p\}.
\]

Fix a point $q\in\Omega_1\setminus\{p\}$.  We claim that there exists a uniform constant $C$ such that
\begin{equation}\label{G(q)leqC}
	G_j(q)\leq C.
\end{equation}

We argue by contradiction and suppose that the non-decreasing sequence $G_j(q) \to\infty$.   Then, the  normalized  functions $v_j=G_j/G_j(q)$ satisfy
\[
 Lv_j=0\quad\text{in }\Omega_j\setminus\{p\},
 \qquad v_j(q)=1.
\]
Then the Harnack inequalities and interior
elliptic estimates give a subsequence converging in
$C^\infty_{\mathrm{loc}}(M\setminus\{p\})$ to a positive solution $v$ of
$$Lv=0,\quad   v(q)=1.$$
Choose a small geodesic
ball $U\Subset\Omega_1$ centered at $p$ with $q\notin\overline U$. Then the  functions
\[
 h_j(x)=\frac{G_j(x)-G_{L,U}(p,x)}{G_j(q)}
\]
are smooth solutions  satisfying $Lh_j=0$ in $U$. 
By the maximum principle, one has $h_j\geq 0$. Using $v_j(q)=1$,  the Harnack inequality yields 
\begin{equation}\label{v_ileqConpartialU}
	C^{-1}\leq v_j\leq C, \quad \mathrm{on} \; \partial U.
\end{equation}
Since $h_j=v_j$ on $\partial U$,  the maximum principle gives
$0\le h_j\leq C$ in $U$. Interior estimates therefore yield a smooth
subsequential limit $h$ in $U$. Since
\[
 v_j=h_j+\frac{G_{L,U}(p,\cdot)}{G_j(q)}\quad\text{on }U\setminus\{p\},
\]
this limit $h$ agrees with $v$ away from $p$ since $G_j(q)\to \infty$. Thus $v$ extends smoothly across $p$ and
\begin{equation}\label{eq:subharmonic}
 \Delta_gv=Vv\ge0\quad\text{on }M.
\end{equation}
Since the function $v_j$ is subharmonic in $\Omega_j\setminus\overline U$ and vanishes on $\partial\Omega_j$, together with \eqref{v_ileqConpartialU}, one has 
\[
 v_j\leq\max_{\partial U}v_j\leq C
 \quad\text{in }\Omega_j\setminus U.
\]
Passing to the limit gives 
$$v\leq  \max_{\partial U}v, \quad \mathrm{in}\;\;M\setminus U.$$
Meanwhile, \eqref{eq:subharmonic} yields $v\leq  \max_{\partial U}v$ in $U$. Consequently,
$v$ attains its global maximum at an interior point of $M$.
The strong maximum principle implies $v\equiv v(q)=1$ which deduces
$$V=Lv=0$$
contradicting $V\not\equiv0.$  Thus, we prove the claim \eqref{G(q)leqC}.

Then, applying Harnack inequalities, interior
estimates, and monotonicity,  we obtain
\[
 G_j\to  G_L
 \quad\text{in }C^\infty_{\mathrm{loc}}(M\setminus\{p\}),
\]
where $G_L$ is finite and positive. The functions $G_j-G_{L,U}$ are smooth, nonnegative, and $L(G_j-G_{L,U})=0$ in $U$, with uniformly bounded boundary
values. Therefore
\[
 G_L=G_{L,U}+h\quad\text{in }U\setminus\{p\}
\]
for a smooth solution $Lh=0$ in $U$. In particular,
$$LG_L=\delta_p$$ distributionally in $M$. Since every Dirichlet Green function is symmetric, its exhaustion limit is symmetric as well.

If $\widetilde G$ is another positive Green function with the same pole
and normalization, then $\widetilde G-G_j$ is smooth and $L(\widetilde G-G_j)=0$
in $\Omega_j$, including at $p$, and is nonnegative on
$\partial\Omega_j$. Thus the maximum principle yields $\widetilde G-G_j\geq 0$ in $\Omega_j$ and then letting $j\to\infty$ deduces $\widetilde G\ge G_L$ which shows that $G_L$ is minimal.

Finally, the maximum principle 
yields 
\[
 G_j\le\max_{\partial U}G_j\le\max_{\partial U}G_L, \quad \mathrm{in}\;\;\Omega_j\setminus \overline{U}
\]
Letting $j\to\infty$, we  prove \eqref{eq:exterior}.
\end{proof}

\begin{lemma}\label{lem:Green-behavior}
Under the same assumption as  in Lemma \ref{lem:green}, let
$\omega_{n-1}$ denote the volume of standard sphere $\ms^{n-1}$ and $r(x)=d_g(p,x)$. Then
\begin{equation}\label{eq:leading}
 G_L(p,x)=\frac{1}{(n-2)\omega_{n-1}}r(x)^{2-n}(1+o(1))
 \qquad\mathrm{as} \;\; x\to p.
\end{equation}
Moreover, in normal coordinates centered at $p$,
\begin{equation}\label{eq:annular-asymptotics}
s^{n-2}G_L\bigl(p,\exp_p(sz)\bigr)
 \longrightarrow\frac{1}{(n-2)\omega_{n-1}}|z|^{2-n}
 \quad\text{as }s\to0
\end{equation}
 in $C^2_{\mathrm{loc}}(\mr^n\setminus\{0\})$.
In particular, there exist positive constants $r_0,C$ such that, for all $0<r(x)<r_0$,
\begin{align*}
 C^{-1}r^{2-n}\le G_L\le Cr^{2-n}; \quad &|\nabla G_L|_g\le Cr^{1-n};\\
 |\nabla\log G_L|_g^2\le Cr^{-2};  \quad & |\nabla^2G_L|_g\le Cr^{-n}.
\end{align*}
\end{lemma}

\begin{proof}
	Following the similar argument as in Lemma \ref{lem:green},
choose a small geodesic ball $U$ centered at $p$.
$G_L-G_{L,U}(p,\cdot)$ extends smoothly across $p$ and satisfies
\begin{equation}\label{HlmabdaleqGlmabda}
 0\le G_L-G_{L,U}(p,\cdot)\leq\max_{\partial U}G_L
 \quad\text{in }U.
\end{equation}
The standard local fundamental-solution expansion for a smooth elliptic
operator with principal part $-\Delta_g$ gives
\[
G_{L,U}(p,x)=\frac{1}{(n-2)\omega_{n-1}}r(x)^{2-n}(1+o(1)).
\]
Together with \eqref{HlmabdaleqGlmabda}, we prove
\eqref{eq:leading}.

Consider the exponential map and define
\[
 v_s(z)=s^{n-2}G_L(p,\exp_p(s z)),
 \qquad g_s=s^{-2}(\exp_p\circ D_s)^*g.
\]
where $D_s(z)=sz$.
On each fixed annulus, the leading expansion gives uniform convergence
of $v_s$ to $\frac{1}{(n-2)\omega_{n-1}}|z|^{2-n}$. Moreover, there holds
\[
 -\Delta_{g_s}v_s
 +s^2V(\exp_p(s z))v_s=0.
\]
The metrics $g_s$ converge smoothly to the Euclidean metric on fixed
annuli, and the rescaled potentials converge smoothly to zero.
Interior elliptic estimates on slightly larger annuli upgrade the
uniform convergence to $C^2$ convergence, proving
\eqref{eq:annular-asymptotics}.  The remaining estimates follow directly.
\end{proof}

\section{Hang-Yang-type inequality}\label{sec:HYinequality}

Let $(M^4,g)$ be a complete four-manifold throughout this section.
If it is closed and $R_g>0$, Hang and Yang
\cite[Proposition 5.1]{HY}
established the following distributional identity
\begin{equation}\label{eq:HY}
	P_g(\log G_{p})
	=16\pi^2\delta_p-Q_g-\frac12\bigl|\Ric_{G_{p}^2g}\bigr|_g^2,
\end{equation}
where $G_{p}$ is the positive Green function, with pole at $p$, of 
the conformal Laplacian operator
\begin{equation}\label{L_g}
	L_gu=-6\Delta_gu+R_gu.
\end{equation}
Integrating \eqref{eq:HY} over the closed $(M,g)$, Hang and Yang \cite{HY} showed that 
\begin{equation}\label{HY-identity}
	\int_M Q_g\,\dmu
	+\frac12\int_{M\setminus\{p\}}\bigl|\Ric_{G_{p}^2g}\bigr|_g^2\,\dmu
	=16\pi^2.
\end{equation}
Motivated by  this identity, we explore the similar behavior on open four-manifolds with positive scalar curvature  by employing Green function obtained  from Lemma \ref{lem:green}.

For any  function $u\in C^4_{loc}$,  define the following Paneitz vector field
\begin{equation}\label{eq:Paneitz-vector}
			\mathcal{J}_g(u)=\nabla\Delta_gu+2\Ric_g^\sharp(\nabla u)-\frac23R_g\nabla u.
\end{equation}
Immediately, using \eqref{Paneitz-operator}, one has 
\begin{equation}\label{divJ=P_g}
	\Div \mathcal{J}_g(u)=P_g u.
\end{equation}

  Let $E_g$ be the trace-free Ricci tensor  and $B_g(p,r)$ be the geodesic ball centered at $p$ with radius $r$ with respect to $g$ throughout  this paper.
\begin{lemma}\label{lem:divK}
	Let $U$ be a connected  open set  of $M^4$. For  a positive function $u\in  C_{loc}^{4}(U) $  and a constant $\varepsilon>0$, define the vector field
	\begin{equation}\label{eq:K-eps}
		\begin{split}
			\mathcal{K}_\varepsilon(u)=&\nabla\left(\frac{|\nabla u|^2-(R_g/6)u^2}{(u+\varepsilon)^2}\right)
			-\frac{2u}{(u+\varepsilon)^2}\Ric_g^\sharp(\nabla u)\\
			&+\frac{2R_gu(u+2\varepsilon)}{3(u+\varepsilon)^3}\nabla u
			+\frac{L_gu}{3(u+\varepsilon)^2}\nabla u.
		\end{split}
	\end{equation}
	Then, one has
	\begin{equation}\label{eq:divK-eps}
		\begin{split}
			\Div 	\mathcal{K}_\varepsilon(u)=\left(1-\theta\right)^2\bigg[Q_g
			&+\frac12\big|E_{u^2g}-4\theta\ud (\log u)\otimes\ud (\log u)+\theta |\nabla \log u|^2g\big|_g^2
			\\
			&+\frac{\theta R_g^2}{18}+\frac{\theta R_g}{3}(1+3\theta)|\nabla\log u|^2
			-\frac{\theta  R_gL_gu}{18u}-\frac{(L_gu)^2}{24u^2}\bigg].
		\end{split}
	\end{equation}
where $\theta=\frac{\varepsilon}{u+\varepsilon}$.
\end{lemma}

\begin{proof}
	For brevity, set $$w=\log u, \quad q=\frac{L_gu}{u}, \quad z=\frac{u}{u+\varepsilon}.$$
	Then,  a direct computation yields 
	\begin{equation}\label{eq:Deltaw}
		\Delta_gw=\frac{R_g-q}{6}-|\nabla w|^2,\qquad R_{u^2g}=qu^{-2}.
	\end{equation}
The conformal transformation formula of Ricci tensor yields
	\begin{equation}\label{conformal-Ricci}
			E_{u^2g}=\Ric_g-2\nabla^2w+2\ud w\otimes\ud w
			-\left(\frac{R_g}{6}+\frac{q}{12}+|\nabla w|^2\right)g.
	\end{equation}
Meanwhile,
\begin{equation}\label{eq:u^2Q_u^2g}
	\begin{split}
			u^4Q_{u^2g}	=&\frac{u^4}{6}\left(\frac{1}{6}L_{u^2g}(R_{u^2g})-3|E_{u^2g}|_{u^2g}^2+\frac{R_{u^2g}^2}{12}\right)\\
		=&\frac{1}{36}uL_g(qu^{-1})-\frac{1}{2}|E_{u^2g}|^2_g+\frac{q^2}{72}\\
		=&-\frac12|E_{u^2g}|_g^2-\frac16\Delta_gq
		+\frac13\langle\nabla w,\nabla q\rangle
		+\frac{R_gq}{18}-\frac{q|\nabla w|^2}{3}-\frac{q^2}{72}.
	\end{split}
\end{equation}
An equivalent expression for \eqref{eq:K-eps} is
	\begin{equation}\label{eq:current-invariant}
			\mathcal{K}_\varepsilon(u)=\left(\frac{R_g}{6}+|\nabla w|^2\right)\nabla(z^2)
		-z^2\left(\mathcal{J}_g(w)+\frac{1}{6}\nabla q-\frac{q}{3}\nabla w\right).
	\end{equation}
	As functions of $w$,
	the derivatives of $z^2$ are
	$$
	(z^2)'=2\theta z^2,
	\qquad
	(z^2)''=(6\theta^2-2\theta)z^2.
	$$
which yields 
	$$
	\nabla(z^2)=(z^2)'\nabla w,
	\qquad
	\Delta_g(z^2)
	=(z^2)''|\nabla w|^2+(z^2)'\Delta_gw.
	$$
Combining with \eqref{conformal-Ricci}, we obtain
	\begin{align*}
		&\nabla\left(\frac{R_g}{6}+|\nabla w|^2\right)
		-\mathcal{J}_g(w)-\frac16\nabla q+\frac q3\nabla w\\
		=&-2E_{u^2g}^\sharp(\nabla w)
		+\left(\frac{R_g}{3}+\frac q6+2|\nabla w|^2\right)\nabla w.
	\end{align*}
Using  \eqref{divJ=P_g}, it follows that
	\begin{align*}
		&\Div	\mathcal{K}_\varepsilon(u)\\
		=&
		\left(\frac{R_g}{6}+|\nabla w|^2\right)\Delta_g(z^2)
		-2E_{u^2g}(\nabla w,\nabla(z^2))\\
		&+\left(\frac{R_g}{3}+\frac q6+2|\nabla w|^2\right)
		\langle\nabla w,\nabla(z^2)\rangle\\
		&+\left(-P_gw-\frac{1}{6}\Delta_gq
		+\frac{1}{3}\langle\nabla w,\nabla q\rangle
		+\frac{q}{3}\Delta_gw\right)z^2\\
		=&	\left(\frac{R_g}{6}+|\nabla w|^2\right)\left((6\theta^2-2\theta)z^2|\nabla w|^2+2\theta z^2(\frac{R_g-q}{6}-|\nabla w|^2)\right)\\
		&-4\theta z^2E_{u^2g}(\nabla w,\nabla w)+2\theta z^2\left(\frac{R_g}{3}+\frac q6+2|\nabla w|^2\right)|\nabla w|^2\\
		&+\left(-P_gw-\frac{1}{6}\Delta_gq
		+\frac{1}{3}\langle\nabla w,\nabla q\rangle
		+\frac{q}{3}\Delta_gw\right)z^2\\
		=&\theta z^2\left((\theta+\frac{1}{3})R_g|\nabla w|^2+6\theta |\nabla w|^4+\frac{R_g(R_g-q)}{18}\right)\\
		&-4\theta z^2E_{u^2g}(\nabla w,\nabla w)+\left(-P_gw-\frac{1}{6}\Delta_gq
		+\frac{1}{3}\langle\nabla w,\nabla q\rangle
		+\frac{q}{3}\Delta_gw\right)z^2\\
	\end{align*}
Using \eqref{conformal-Q} and \eqref{eq:u^2Q_u^2g}, we have
	$$
	\begin{aligned}
		&-P_gw-\frac16\Delta_gq
		+\frac13\langle\nabla w,\nabla q\rangle
		+\frac q3\Delta_gw\\
		=&Q_g+\frac12|E_{u^2g}|_g^2
		-\frac{R_gq}{18}
		+\frac{q|\nabla w|^2}{3}
		+\frac{q^2}{72}
		+\frac q3\left(\frac{R_g-q}{6}-|\nabla w|^2\right)\\
		=&Q_g+\frac12|E_{u^2g}|_g^2-\frac{q^2}{24}.
	\end{aligned}
	$$
	Combining these identities,  one has 
	\begin{align*}
		&\Div 	\mathcal{K}_\varepsilon(u)\\
		=&
		z^2\Bigg[
		Q_g+\frac12|E_{u^2g}|_g^2
		-4\theta E_{u^2g}(\nabla w,\nabla w)
		+6\theta^2|\nabla w|^4\\
		&\qquad\quad
		+\theta\left(\theta+\frac13\right)R_g|\nabla w|^2
		+\frac{\theta}{18}R_g(R_g-q)
		-\frac{q^2}{24}
		\Bigg].
	\end{align*}
	Finally, completing the square gives \eqref{eq:divK-eps}.
\end{proof}

In this section, from now on, we deal with the case 
\begin{equation}\label{scalar-condition}
 Q_g\geq 0, \quad  R_g>0.
\end{equation}
All gradients, divergences, and tensor norms in this section are taken
with respect to $g$, unless another metric is indicated.

Fix  a point $p\in M$. By  Lemma \ref{lem:green}, the operator $L_g$ \eqref{L_g} has a   positive  minimal Green function $G_p$ with pole $p$. We normalize it  by $$L_g G_p=\delta_p$$ and set $G=G_p$ if no confusion is caused throughout this paper.

In the estimates below, $K$ is a fixed compact neighborhood of $p$. The constant $C$ may be different from line to line and   may depend on $g$, $p$, $K$, the fixed inner cutoffs.
\begin{lemma}\label{lem:weighted}
Under  \eqref{scalar-condition}, one has 
\begin{align}
 \int_{M\setminus K}\left(|\nabla G|^2
       +R_gG^2\right)\dmu&\leq C,
       \label{eq:first}\\
 \int_{M\setminus K}\left[R_g|\nabla G|^2
       +(R_g^2+|\Ric_g|^2+Q_g)G^2\right]\dmu&\leq C,
       \label{eq:curvature}\\
      \int_{M\setminus K}G^2|\nabla R_g|\dmu
     &\leq C,  \label{eq:Rgradient}\\
     	\int_{M\setminus K}\left(
     	|\nabla^2G|^2+
     	\frac{|\nabla G|^4}{G^2}\right)\dmu &\leq C.\label{eq:hessian}
\end{align}
\end{lemma}

\begin{proof}
Let $\Omega$ range over a smooth relatively compact exhaustion of $M$
containing $K$, and let $G_{\Omega}$ be the Dirichlet Green
function of $L_g$ on $\Omega$. 
For a fixed smaller  compact neighborhood $p\in K_0\Subset K$, choose a smooth function
$0\leq\chi\leq1$ satisfying
\begin{equation}\label{chi-test}
	\chi\equiv 0,\quad \mathrm{in}\;\;  K_0; \quad \chi \equiv1,\quad  \mathrm{in}\;\;  M\setminus K.
\end{equation}
Thus $\nabla\chi$ is supported in a compact set contained in $K\setminus K_0$. Testing $L_g G_{\Omega}=0$ away from the pole with
$\chi^2 G_{\Omega}$ gives
\begin{equation}\label{eq:caccioppoli}
 6\int_\Omega|\nabla(\chi G_{\Omega})|^2\dmu
 +\int_\Omega R_g\chi^2G_{\Omega}^2\dmu
 =6\int_\Omega G_{\Omega}^2|\nabla\chi|^2\dmu.
\end{equation}
Using the fact $0\leq G_{\Omega}\leq G$, the right-hand side is controlled by
$$\int_\Omega G_{\Omega}^2|\nabla\chi|^2\dmu\leq \int_{K\setminus K_0}G^2|\nabla\chi|^2\dmu\leq C$$
which is independent of $\Omega$.
 Exhausting $M$ and applying Fatou's lemma, we prove \eqref{eq:first}.

In $\Omega\setminus\{p\}$, a direct computation yields 
\begin{equation}\label{eq:div_gR_gG^2}
	\begin{split}
		&\Div\bigl(R_g\nabla(G_{\Omega}^2)-G_{\Omega}^2\nabla R_g\bigr)\\
		=&2R_g|\nabla G_{\Omega}|^2+\left(\frac1{12}R_g^2+3|E_g|^2+6Q_g\right)G_{\Omega}^2.
	\end{split}
\end{equation}

Since $G_\Omega=0$ on $\partial \Omega$, one has 
$$R_g\nabla(G_{\Omega}^2)-G_{\Omega}^2\nabla R_g=2G_\Omega R_g\nabla(G_{\Omega})-G_{\Omega}^2\nabla R_g=0, \quad \mathrm{on}\;\; \partial \Omega.$$
Then,  multiplying \eqref{eq:div_gR_gG^2} by $\chi^2$ and integrating by parts,  we obtain
\begin{align*}
 &\int_\Omega\chi^2\left[
 2R_g|\nabla G_{\Omega}|^2+
 \left(\frac1{12}R_g^2
       +3|E_g|^2+6Q_g\right)G_{\Omega}^2\right]\dmu\\
 =&-4\int_\Omega\chi R_gG_{\Omega}\ip{\nabla\chi}{\nabla G_{\Omega}}\dmu
       +2\int_\Omega\chi G_{\Omega}^2\ip{\nabla\chi}{\nabla R_g}\dmu.
\end{align*}
We absorb the first term by using Young's inequality
\[
 4\chi R_gG_{\Omega}|\nabla\chi|\cdot|\nabla G_{\Omega}|
 \leq\chi^2R_g|\nabla G_{\Omega}|^2+4R_gG_{\Omega}^2|\nabla\chi|^2.
\]
Meanwhile, 
$$4\int_\Omega R_gG_{\Omega}^2|\nabla\chi|^2\dmu\leq4\int_{K\setminus K_0}R_gG^2|\nabla\chi|^2\dmu\leq C.$$
Similarly, there holds
\[2\int_\Omega\chi G_{\Omega}^2\ip{\nabla\chi}{\nabla R_g}\dmu\leq 2 \int_{K\setminus K_0}|\nabla\chi|\cdot|\nabla R_g|G^2\dmu\leq C.\]
Passing to the
exhaustion limit, we prove
\eqref{eq:curvature}.

A direct computation yields
\[
 \Delta_g\log R_g
 =\frac{R_g}{4}-
   \frac{3|E_g|^2+6Q_g}{R_g}-|\nabla\log R_g|^2\leq \frac{R_g}{4} -|\nabla\log R_g|^2.
\]
Multiplying the above inequality by a cutoff function $\zeta^2$,  integrating by parts, and using Young's inequality, one has
\begin{equation*}
 \frac12\int_M\zeta^2|\nabla\log R_g|^2\dmu
 \leq2\int_M|\nabla\zeta|^2\dmu
       +\frac14\int_MR_g\zeta^2\dmu.
\end{equation*}
 Apply this inequality
to $\zeta=\chi G_{\Omega}$, extended by zero outside $\Omega$.  Applying the estimate
\eqref{eq:caccioppoli} and  exhaustion,  Fatou's lemma yields  
\begin{equation*}
	 \int_{M\setminus K}G^2
	|\nabla\log R_g|^2\dmu\leq C.
\end{equation*}
Combining with  the estimates \eqref{eq:curvature}, H\"older's inequality gives 
\begin{align*}
 &\int_{M\setminus K}G^2|\nabla R_g|\dmu\\
 \leq &
 \left(\int_{M\setminus K}G^2
             |\nabla\log R_g|^2\dmu\right)^{\frac{1}{2}}\cdot
 \left(\int_{M\setminus K}R_g^2G^2\dmu\right)^{\frac{1}{2}}\leq C
\end{align*}
which proves \eqref{eq:Rgradient}. 

Thus, the estimates  \eqref{eq:first}, \eqref{eq:curvature} and   \eqref{eq:Rgradient} still hold on $M\setminus K_0$ by the same argument.

With the help of the estimate  \eqref{eq:curvature} and 
$\Delta_gG=R_gG/6$ on $M\setminus\{p\}$, we obtain
\begin{equation*}
 \int_{M\setminus K}\left((\Delta_gG)^2+G^2|\Ric_g|^2\right)\dmu
 \leq C.
\end{equation*}

 Completeness of metric provides
compactly supported Lipschitz functions $0\leq\eta_j\leq1$ such that
\begin{equation}\label{cut-eta_j}
\eta_j\equiv 1,\quad \mathrm{in}\;\;  B_g(p,j); \quad \eta_j\equiv0,\quad  \mathrm{in}\;\;  M\setminus B_g(p,2j);  \quad |\nabla\eta_j|\leq \frac{C}{j}.
\end{equation}
  Then,  using \eqref{eq:first}, one has
\begin{align*}
	&\int_M|\nabla(\chi\eta_j) |^2|\nabla G|^2\dmu\\
	\leq &2\int_M|\nabla\chi|^2\eta_j^2|\nabla G|^2\dmu+2\int_M\chi^2|\nabla \eta_j|^2|\nabla G|^2\dmu\\
	\leq &C\int_{M\setminus K_0}|\nabla G|^2\dmu \leq C.
\end{align*}
Combining these two estimates, we obtain
\begin{equation}\label{eq:cutoff-energy}
 \int_M\chi^2\eta_j^2\left((\Delta_gG)^2+G^2|\Ric_g|^2\right)\dmu
 +\int_M|\nabla(\chi\eta_j)|^2|\nabla G|^2\dmu\leq C
\end{equation}
uniformly in $j$. 

For brevity, set
\[
 A_1=\int_M\chi^2\eta_j^2|\nabla^2G|^2\dmu,
 \qquad A_2=\int_M\chi^2\eta_j^2\frac{|\nabla G|^4}{G^2}\dmu.
\]
With the help of Bochner's formula and integration by parts, one has
\begin{align*}
 A_1=&\int_M\chi^2\eta_j^2(\Delta_gG)^2\dmu
      -\int_M\chi^2\eta_j^2\Ric_g(\nabla G,\nabla G)\dmu\\
 &+2\int_M\chi\eta_j\Delta_gG\langle \nabla(\chi\eta_j), \nabla G \rangle \dmu
      -2\int_M\chi\eta_j\nabla^2G(\nabla(\chi\eta_j),\nabla G)\dmu.
\end{align*}
Using \eqref{eq:cutoff-energy}, H\"older's inequality and Young's inequality, we  have
\[
\left|\int_M\chi^2\eta_j^2\Ric_g(\nabla G,\nabla G)\dmu\right|\leq  \left(\int_M\chi^2\eta_j^2G^2|\Ric_g|^2\dmu\right)^{\frac{1}{2}}\cdot A_2^{\frac{1}{2}}
 \leq C A_2^{\frac{1}{2}},
\]
\[2\int_M\chi\eta_j\Delta_gG\langle\nabla(\chi\eta_j), \nabla G \rangle \dmu\leq \int_M\chi^2\eta_j^2(\Delta_gG)^2\dmu+\int_M|\nabla(\chi\eta_j)|^2|\nabla G|^2\dmu\leq C,\]
\[-2\int_M\chi\eta_j\nabla^2G(\nabla(\chi\eta_j),\nabla G)\dmu\leq \frac{1}{2}A_1+2\int_M|\nabla(\chi\eta_j)|^2|\nabla G|^2\dmu\leq\frac{1}{2}A_1+C. \]
Thus, we conclude that
\begin{equation}\label{eq:H-bound}
 A_1\leq C A_2^{\frac{1}{2}}+C.
\end{equation}
On the other hand, on $M\setminus\{p\}$, a direct computation yields
\[
 \Div\left(\frac{|\nabla G|^2}{G}\nabla G\right)
 =\frac{2\nabla^2G(\nabla G,\nabla G)
         +(\Delta_gG)|\nabla G|^2}{G}
   -\frac{|\nabla G|^4}{G^2}.
\]
Multiplying by $\chi^2\eta_j^2$ and integrating gives
\begin{align*}
 A_2\leq{}&2A_1^{\frac{1}{2}}A_2^{\frac{1}{2}}
 +\left(\int_M\chi^2\eta_j^2(\Delta_gG)^2\dmu\right)^{\frac{1}{2}}A_2^{\frac{1}{2}}\\
 &+2\left(\int_M|\nabla(\chi\eta_j)|^2|\nabla G|^2\dmu\right)^{\frac{1}{2}}A_2^{\frac{1}{2}}\\
 \leq&(2A_1^{\frac{1}{2}}+C)A_2^{\frac{1}{2}}.
\end{align*}
Thus, squaring both sides, dividing both sides by $A_2$, and then using \eqref{eq:H-bound}, we obtain $$A_2\leq CA_1+C\leq CA_2^{\frac{1}{2}}+C.$$ Thus, both
$A_1$ and $A_2$ are bounded  independently of $j$. Letting $j\to\infty$
and applying Fatou's lemma proves \eqref{eq:hessian}.
\end{proof}

Now,  we  are ready to establish the following Hang-Yang-type inequality.

\begin{theorem}\label{thm:CVandHY-for-Q}
Let $(M^4,g)$ be an open four-manifold with  $Q_g\geq 0$ and
	$R_g>0$. Then 
	\begin{equation}\label{eq:final-g}
		\int_MQ_g\dmu
		+\frac12\int_{M\setminus\{p\}}|\Ric_{G_p^2g}|_g^2\dmu
		\leq16\pi^2.
	\end{equation}
For any constant $\eps>0$, there holds
\begin{equation}\label{intR^2}
	\int_{M\setminus\{p\}}R_g^2\frac{\eps \cdot G_p^2}{(G_p+\eps)^3}\dmu\leq 288 \pi^2.
\end{equation}
\end{theorem}

\begin{proof}
 On $M\setminus\{p\}$ and $\eps>0$, consider the vector field   $\mathcal{K}_\eps(G)$ defined in \eqref{eq:K-eps}.  Let $K$ be  a compact neighborhood $K$ of $p$. In $M\setminus K$, using \eqref{eq:K-eps},  one has 
 \begin{align*}
 	|\mathcal{K}_\eps(G)|\leq C_\eps\bigg(&
 	|\nabla^2G|\cdot |\nabla G|+\frac{|\nabla G|^3}{G}
 	+G^2|\nabla R_g|\\
 	&+G|\Ric_g|\cdot |\nabla G|+R_gG|\nabla G|\bigg).
 \end{align*}
 Each term on the right is integrable on $M\setminus K$ by
 Lemma \ref{lem:weighted} and H\"older's
 inequality. Thus, one has
 \begin{equation}\label{exterior-Jr}
 	\int_{M\setminus K}|\mathcal{K}_\eps(G)|\dmu\leq C(\varepsilon).
 \end{equation}
Lemma \ref{lem:Green-behavior} shows that 
\begin{equation}\label{eq:pole}
s^2G \bigl(\exp_p(s z)\bigr)
 \longrightarrow\frac1{24\pi^2|z|^2}
 \quad\text{in }C^2_{\mathrm{loc}}(\mr^4\setminus\{0\}),
 \quad \mathrm{as}\; \; s\to 0.
\end{equation}
Let $r=d_g(p,x)$. 
Using \eqref{eq:K-eps} and \eqref{eq:pole}, near the pole $p$, one has
\begin{equation}\label{eq:near-pole}
 \mathcal{K}_\eps(G)=-8r^{-3}\partial_r+o(r^{-3}).
\end{equation}
Thus $\mathcal{K}_\eps(G)$ is locally integrable near $p$, and combining with \eqref{exterior-Jr}, one has
\begin{equation}\label{eq:KinL1}
	\mathcal{K}_\eps(G)\in L^1(M,g)
\end{equation}
for each fixed  $\eps>0$.

For brevity, set the notations
\begin{equation*}
	w =\log G,\qquad
	\theta=\frac{\eps}{G+\eps}.
\end{equation*}
Using Lemma \ref{lem:divK},  on $M\setminus\{p\}$, we have 
\begin{equation}\label{eq:divK}
	\begin{split}
		\Div \mathcal{K}_\eps(G)=\left(1-\theta\right)^2\bigg[&Q_g+
		\frac12\left|E_{G^2g}-4\theta\ud w\otimes\ud w+\theta |\nabla w|^2 g\right|_g^2\\
		&+\frac{\theta R_g}{3}\left(\frac{R_g}{6}+(1+3\theta)|\nabla w|^2\right)\bigg].
	\end{split}
\end{equation}
Moreover, using \eqref{eq:near-pole},  one has
\begin{equation}\label{eq:pole-flux}
		\lim_{r\to 0}\int_{\partial B_g(p,r)}\langle\mathcal{K}_\eps(G), \partial r\rangle\ud\sigma_g=-16\pi^2.
\end{equation}
For all nonnegative $\eta\in C^\infty_c(M)$, we obtain the distributional inequality
\begin{equation}\label{eq:distribution}
 \begin{split}
 &\int_{M\setminus\{p\}}\eta\cdot \left(1-\theta\right)^2\left[Q_g+\frac{\theta R_g^2}{18}+
  \frac12\left|E_{G^2 g}-4\theta \ud w\otimes\ud w+\theta |\nabla w|^2 g\right|_g^2\right]\dmu\\
\leq & -\int_M\langle \mathcal{K}_\eps(G),\nabla\eta\rangle\dmu+16\pi^2\eta(p).
 \end{split}
\end{equation}
Testing  \eqref{eq:distribution} by  the  cutoff functions $\eta_j$ from \eqref{cut-eta_j} and using \eqref{eq:KinL1},  one has 
\[
 \left|\int_M\langle \mathcal{K}_\eps(G), \nabla\eta_j\rangle \dmu\right|
 \leq\frac Cj\int_M|\mathcal{K}_\eps(G)|\dmu\longrightarrow0,\quad \mathrm{as}\quad  j\to\infty. 
\]
Fatou's lemma on the left and the dominated convergence   on the right  in \eqref{eq:distribution} yield
\begin{equation}\label{eq:finiteE_G}
 \begin{split}
 \int_{M\setminus\{p\}}\left(1-\theta\right)^2\left[Q_g+ \frac{\theta R_g^2}{18}+
   \frac12\left|E_{G^2g}-4\theta\ud w\otimes\ud w +\theta |\nabla w|^2 g\right|_g^2\right]\dmu\leq 16\pi^2.
 \end{split}
\end{equation}
The estimate \eqref{intR^2} follows immediately. For any $x\neq p$, $G(x)>0$ and hence  $\theta(x)\to0$ as
$\eps\to0$. Applying  Fatou's lemma again,  we prove \eqref{eq:final-g}.
\end{proof}

\section{Finiteness of fundamental group }\label{sec:pi1}

The  Cohn-Vossen-type inequality \eqref{eq:final-g} implies the finiteness of $\pi_1(M)$ when $Q_g\not \equiv 0.$

\begin{proposition}\label{prop:finitepi}
Let $(M^4,g)$ be an open four-manifold with
\[
  R_g>0,\qquad Q_g\geq0,\qquad Q_g\not\equiv0.
\]
Then $M$ has one end and  the fundamental group $\pi_1(M)$ is finite.
\end{proposition}

\begin{proof}
	An open  Riemannian manifold has at least one
	end. If $M$ had more than one end,
	Theorem \ref{thm:splitting} would force $Q_g\equiv0$, a
	contradiction. Hence $M$ has exactly one end. 
	
Let $\pi:(\widetilde M,\tilde g)\to(M,g)$ be the universal
cover, with $\tilde g=\pi^*g$. 
First, by the properties  of the universal cover,  $\tilde g$ is complete and 
$$R_{\tilde g}> 0, \quad Q_{\tilde g}\geq0,\quad  Q_{\tilde g}\not\equiv 0.$$
 Theorem \ref{thm:CVandHY-for-Q} yields 
 \begin{equation}\label{eq:cover-Q-bound}
 	\int_{\widetilde M}Q_{\tilde g}\,\ud\mu_{\tilde g}
 	\leq16\pi^2.
 \end{equation}

Since $Q_g\not\equiv 0$, there exists a point $x_0\in M$ such that $Q_g(x_0)>0$ and then 
choose a small geodesic  ball  $ B_g(x_0,r_0)\subset M$  such that  $\pi^{-1}(B_g(x_0,r_0))$ is the disjoint union of isometric lifts of $B_g(x_0,r_0)$.  Meanwhile, in $B_{g}(x_0,r_0)$, one has
$$\int_{B_g(x_0,r_0)}Q_g\dmu>0.$$
Choosing  $m$ disjoint lifts of $B_g(x_0,r_0)$, one has 
\[
  16\pi^2\geq \int_{\widetilde M}Q_{\tilde g}\,\ud\mu_{\tilde g}
  \geq m \int_{B_g(x_0,r_0)} Q_g\dmu
\]
which yields   $\pi_1(M)$ must be finite.  
\end{proof}

Now, we are going to deal with the case of vanishing $Q$-curvature. 
 First, we recall a corollary of Brinkmann's theorem \cite{Brin}.
For the reader's convenience, we include the proof here.
\begin{lemma}\label{lem:Brinkman}
	Let $(M,g)$ be a four-manifold. Suppose that $g$ is  Ricci-flat and $v$ is a positive smooth function in a connected neighborhood  $U$ of some point. If  the conformal metric $v^{-2}g$ is also
	Ricci-flat in $U$, then either $v$ is constant or $g$ is flat in $U$.
\end{lemma}

\begin{proof}
The conformal Ricci formula yields
$$\Ric_{v^{-2}g}=\Ric_g+2\frac{\nabla^2 v}{v}+\left(\frac{\Delta_g v}{v}-3\frac{|\nabla v|^2}{v^2}\right)g.$$
Using Ricci-flat assumption and taking trace, one has
\begin{equation}\label{trace-conformal}
	\frac{\Delta_g v}{v}=2\frac{|\nabla v|^2}{v^2}.
\end{equation}
Inserting it back,  a direct computation yields
	\begin{equation}\label{eq:concircular}
		\nabla^2v=\frac{|\nabla v|^2}{2v}g.
	\end{equation}
On  the one hand, one has
$$\mathrm{div}_g(\nabla^2 v)=\mathrm{div}_g(\frac{|\nabla v|^2}{2v} g)=\nabla \frac{|\nabla v|^2}{2v} .$$
On the other hand, using \eqref{trace-conformal} and the  Ricci identity, one has 
	$$\mathrm{div}_g	(\nabla^2v)=\nabla (\Delta_g v)+\Ric_g^{\sharp}(\nabla v)=4\nabla \frac{|\nabla v|^2}{2v} .$$
Combining these two identities, we have  $\nabla \frac{|\nabla v|^2}{2v} =0$ and then $\frac{|\nabla v|^2}{2v} $ is a nonnegative constant. If
$\frac{|\nabla v|^2}{2v} =0$, then $v$  is a constant.
Otherwise, using  \eqref{eq:concircular}, one has
$\nabla^3 v=0.$
Then, for all vector fields $Y,Z$, one has
	\[
	\mathrm{Rm}_g(Y,Z)\nabla v=\nabla_Y\nabla_Z\nabla v-\nabla_Z\nabla_Y\nabla v-\nabla_{[Y,Z]}\nabla v=0.
	\]
In $U$, all curvature components involving the unit vector
	parallel to \(\nabla v\) vanish. The remaining components form an
	algebraic curvature tensor on its three-dimensional orthogonal
	complement, with zero Ricci contraction. Such a tensor is zero in
	dimension three. Hence $g$ is flat.
\end{proof}

\begin{lemma}\label{lem:lcf}
	Let  $(M^4,g)$  be a complete four-manifold with 
	$$Q_g\equiv0, \quad  R_g>0.$$
	If the fundamental group $\pi_1(M)$ is infinite, then the Weyl tensor $W_g\equiv 0.$
\end{lemma}

\begin{proof}
	Let  $(\widetilde M, \tilde g) $ be the universal  cover. Since $\pi_1(M)$ is infinite, $\widetilde M $ must be noncompact with $R_{\tilde g}>0$.  For brevity,  let $G=G_{L_{\tilde g}}(p,\cdot)$ be the minimal positive Green function of $L_{\tilde g}$ with  the pole $p$  using Lemma \ref{lem:green}.
	
	  Let $\Gamma$ be the group of  deck transformation which is isomorphic to $\pi_1(M)$.
	Choose a sequence of distinct nonidentity elements $\gamma_j\in\Gamma$, and define
	\begin{equation}\label{eq:normalized}
		u_j(x)=\frac{G(\gamma_j x)}{G(\gamma_j p)}.
	\end{equation}
	The points $\gamma_j^{-1}p$, which are the poles of $u_j$,
	leave every compact subset as $j\to\infty$. On each fixed relatively compact
	set, $u_j$ is eventually a positive solution to
	 $$L_{\tilde g}u_j=0,\quad u_j(p)=1.$$  Up to a subsequence, the local Harnack inequalities and interior
	elliptic estimates deduce
	\begin{equation}\label{eq:limit-u}
		u_j\longrightarrow u \text{ in } C^\infty_{loc}(\widetilde M),
		\qquad u>0,\quad u(p)=1,\quad L_{\tilde g}u=0.
	\end{equation}
	For any compact set $K\subset \widetilde M $ and   all sufficiently large $j$, isometry invariance yields 
	\[
	\int_K|\Ric_{u_j^2\tilde g}|_{\tilde g}^2\ud\mu_{\tilde g}
	=\int_{\gamma_jK}|\Ric_{G^2\tilde g}|_{\tilde g}^2\ud\mu_{\tilde g}.
	\]
	Properness of the deck action implies that \(\gamma_jK\) leaves
	every fixed compact set. Then the  right-hand side  tends to zero by
	\eqref{eq:final-g}. Therefore, combining with \eqref{eq:limit-u},
$h=u^2\tilde g$ is Ricci-flat.

Applying \eqref{intR^2} with $\eps=G(\gamma_j p)$ and invariant of $\tilde g$ under $\gamma_j$, one has
\begin{align*}
	\int_{\widetilde M} R_{\tilde g}^2\frac{u_j^2}{(u_j+1)^3}\ud\mu_{\tilde g}
	=&\int_{\widetilde M} R_{\tilde g}^2(x)\frac{G(\gamma_j p)G^2(\gamma_jx)}{(G(\gamma_jx)+G(\gamma_j p))^3}\ud\mu_{\tilde g}\\
	=&\int_{\widetilde M} R_{\tilde g}^2\frac{G(\gamma_j p)G^2}{(G+G(\gamma_j p))^3}\ud\mu_{\tilde g}\leq288\pi^2.
\end{align*}
Combining with	\eqref{eq:limit-u}, the Fatou's lemma yields 
	\begin{equation}\label{eq:u-weighted}
		\int_{\widetilde M} R_{\tilde g}^2\frac{u^2}{(1+u)^3}\ud\mu_{\tilde g}\leq288\pi^2.
	\end{equation}

	Now, suppose  that $h$ is not flat.
	  
	    If $v$ is  another  limit obtained from another subsequence of  \eqref{eq:normalized}, then  $v^2\tilde g$ is also  Ricci-flat. Since 
	  $  v^2\tilde g=(v/u)^2h$,   Lemma \ref{lem:Brinkman}, the nonflatness of $h$ and  $v(p)=u(p)=1$ imply $v\equiv u$. Thus, the normalized Green functions
	  converge to the same limit $u$ along every such sequence.

	   For each $\gamma\in\Gamma$, the metric
	  \[
	  \gamma^*h=(u\circ\gamma)^2\tilde g
	  =\left(\frac{u\circ\gamma}{u}\right)^2h
	  \]
	  is also Ricci-flat.  Lemma \ref{lem:Brinkman} implies that
	\begin{equation}\label{eq:character}
		u\circ\gamma=c_\gamma u
	\end{equation}
for some positive constant $c_\gamma$.
	The map $\gamma\mapsto c_\gamma$  is a multiplicative character for $\gamma\in \Gamma$.
	
	If $c_\gamma\equiv 1$,
	the strictly positive density $R_{\tilde g}^2\frac{u^2}{(1+u)^3}\,d\mu_{\tilde g}$ is
	$\Gamma$-invariant. Choose a nonempty relatively compact
	open set $U$ whose distinct deck translates are disjoint.
	Then
	\[
	\int_{\widetilde M} R_{\tilde g}^2\frac{u^2}{(1+u)^3}\ud\mu_{\tilde g}
	\geq\sum_{\gamma\in\Gamma}
	\int_{\gamma U}R_{\tilde g}^2\frac{u^2}{(1+u)^3}\ud\mu_{\tilde g}
	=\sum_{\gamma\in\Gamma}\int_U R_{\tilde g}^2\frac{u^2}{(1+u)^3}\ud\mu_{\tilde g}
	=+\infty,
	\]
	contrary to \eqref{eq:u-weighted}. Thus there exists an element 
	$\gamma_0\in\Gamma$  with $c_{\gamma_0}\not=1$ and  then such $\gamma_0$ has infinite order.

	For each non-zero integer $m$, set
	\[
	b_m=G_{L_{\tilde g}}(p,\gamma_0^mp).
	\]
	Symmetry and isometry invariance 
	give $b_m=b_{-m}$. Applying the preceding convergence to
	the sequence $\{\gamma_0^m\}$, and evaluating at $\gamma_0 p$, yields
	\[
	\frac{b_{m+1}}{b_m}
	\to  u(\gamma_0 p)=c_{\gamma_0} u(p)=c_{\gamma_0}
		\qquad(m\to+\infty).
	\]
	Applying it to \(\gamma_0^{-m}\), and evaluating at
	\(\gamma_0^{-1}p\), yields
	\[
	\frac{b_{-m-1}}{b_{-m}}
	\to u(\gamma_0^{-1}p)=c_{\gamma_0}^{-1}.
	\]
	The two ratios are equal for every $m\geq1$. Hence
	$c_{\gamma_0}=c_{\gamma_0}^{-1}$ which contradicts
$c_{\gamma_0}\not=1$.

	It follows that $h$ is flat and in particular, one has $W_h\equiv0$. 
	 Since $h=u^2\tilde g$,  then the conformal covariant property of Weyl tensor  yields 
	 $$W_{\tilde g}\equiv 0.$$
	 Thus, the Weyl tensor of $g$ must also vanish.
\end{proof}
\begin{remark}
	In the closed  oriented case, Gursky showed $W_g\equiv 0$ under the additional assumption $H^1(M,\mr)\not=0$;  see \cite[Corollary F  (ii)]{Gursky Ann}.
\end{remark}

\begin{proposition}\label{prop:finitepiforQ=0}
	Let $(M,g)$ be an open four-manifold with  $R_g>0$ and $Q_g \equiv0$.
	Then  $\pi_1(M)$ is finite.
\end{proposition}
\begin{proof}
	We argue by contradiction. If $\pi_1(M)$ is infinite, by Lemma \ref{lem:lcf}, $(M,g)$ is locally conformally flat. Its universal cover is complete, simply connected, and noncompact. Proposition \ref{prop:cylinder} therefore shows  that the universal cover is isometric to $(\mr\times\ms^3,\ud t^2\oplus a^2g_{\ms^3})$ for some constant $a>0$.
	Each isometry of $\widetilde M$ preserves the one-dimensional
	kernel of its Ricci tensor and its orthogonal complement.
	Every deck transformation therefore has the form
	\[
	\gamma(t,\theta)
	=(\epsilon_\gamma t+b_\gamma,A_\gamma\theta),
	\qquad \epsilon_\gamma\in\{1,-1\},\quad
	A_\gamma\in O(4).
	\]
	The subgroup
	$\Gamma^+=\{\gamma\in\Gamma:\epsilon_\gamma=1\}$
	has index at most two, so it is infinite.
	If $b_\gamma=0$ for all $\gamma\in\Gamma^+$, then
	$\Gamma^+$ preserves the compact set
	$\{0\}\times\ms^3$, contradicting proper
	discontinuity. Choose $\gamma\in\Gamma^+$ with
	$b_\gamma\ne0$. The quotient $\widetilde M/\langle\gamma\rangle$
	is compact since  a closed slab of width $|b_\gamma|$
	maps onto it. The natural map
	\[
	\widetilde M/\langle\gamma\rangle\longrightarrow \widetilde M/\Gamma=M
	\]
	is continuous and surjective. Hence $M$ is compact,
	contrary to our assumption. 
\end{proof}

\section{Modified Green function}\label{sec:modified-Green}

Throughout this section,  let $(M,g)$ be an open   four-manifold with $$Q_g\geq0, \quad R_g>0.$$ 
\begin{lemma}\label{lem:u-modifiedG}
Consider a non-zero and nonnegative function $f\in C_c^{2,\alpha_0}(M)$ for some $0<\alpha_0<1$.  Then, there exists a minimal positive solution $u$ to  $L_gu=f.$
	 Moreover, $$u\leq\sup_M(f/R_g)$$ and
	\begin{equation}\label{eq:potential-energy}
		\int_M\bigl(6|\nabla u|^2+R_gu^2\bigr)\ud\mu_g
		\leq\int_M\frac{f^2}{R_g}\ud\mu_g.
	\end{equation}
	There is a constant $C_f$ depending on $f$  such that
	\begin{equation}\label{eq:potential-curvature}
		\int_M\left[R_g|\nabla u|^2+
		\left(R_g^2+|\Ric_g|^2+Q_g\right)u^2\right]\ud\mu_g\leq C_f
	\end{equation}
	\begin{equation}\label{eq:potential-second}
		\int_M\left(u^2|\nabla R_g|+u^2|\nabla\log R_g|^2+|\nabla^2u|^2+
		\frac{|\nabla u|^4}{u^2}\right)\ud\mu_g\leq C_f
	\end{equation}
\end{lemma}

\begin{proof}
	Let $\Omega\Subset M$ be a smooth connected domain containing
	$\mathrm{supp} f$. There exists a solution $v_\Omega$  to    
	\begin{equation}\label{eq:v-equation}
		L_g v_\Omega=f, \; \mathrm{in} \;\;\Omega; \quad v_\Omega=0, \; \mathrm{on}\;\; \partial \Omega
	\end{equation}
	The maximum
	principle yields  $v_\Omega>0$ in $\Omega$. The finite constant
	$\sup_M(f/R_g)$ is a supersolution, so
	\begin{equation}\label{v-upper-bound}
		v_\Omega\leq\sup_M(f/R_g)
	\end{equation}
	independently of $\Omega$.
	Multiplying \eqref{eq:v-equation} by $v_\Omega$,  integrating by parts and using H\"older's inequality, one has 
	\begin{equation*}
		\int_\Omega\bigl(6|\nabla v_\Omega|^2+R_gv_\Omega^2\bigr)\ud\mu_g
		=\int_\Omega fv_\Omega\ud\mu_g
		\leq\left(\int_M\frac{f^2}{R_g}\ud\mu_g\right)^{1/2}
		\left(\int_\Omega R_gv_\Omega^2\ud\mu_g\right)^{1/2}
	\end{equation*}
	which deduces $$\int_\Omega R_g v_\Omega^2\dmu\leq \int_M\frac{f^2}{R_g}\ud\mu_g$$ and then
	\begin{equation}\label{eq:nablav2}
		\int_\Omega\bigl(6|\nabla v_\Omega|^2+R_gv_\Omega^2\bigr)\ud\mu_g\leq \int_M\frac{f^2}{R_g}\ud\mu_g.
	\end{equation}
	
	The Dirichlet solutions  to \eqref{eq:v-equation} increase under domain exhaustion. Their
	uniform upper bound  \eqref{v-upper-bound} and local elliptic estimates give a positive
	limit $u$ to the following equation
	\begin{equation*}
		L_g u =f, \quad \mathrm{in}\;\; M.
	\end{equation*}
	Comparison with any other positive solution proves minimality of $u$.
	Using  Fatou's lemma in \eqref{eq:nablav2} as $\Omega$ exhausts $M$, we prove   \eqref{eq:potential-energy}.
	
A direct computation yields
$$\Div(R_g\nabla (v_\Omega)^2-v_\Omega^2\nabla R_g)=2R_g|\nabla v_\Omega|^2+\left(\frac{R_g^2}{12}+3|E_g|^2+6Q_g\right)v_\Omega^2-\frac{1}{3}R_gfv_\Omega.$$
Integrating the above identity over $\Omega$, we obtain
\begin{align*}
&	\int_\Omega\left(2R_g|\nabla v_\Omega|^2+\left(\frac{R_g^2}{12}+3|E_g|^2+6Q_g\right)v_\Omega^2\right)\ud\mu_g\\
=&\frac{1}{3}\int_\Omega R_gfv_\Omega\dmu\leq \frac{1}{3}\sup_M(f/R_g)\int_MfR_g\dmu\leq C_f.
\end{align*}
The exhaustion and Fatou's lemma deduce \eqref{eq:potential-curvature}.

Imitating the proof of Lemma \ref{lem:weighted} and using the fact $f$ has compact support,  it is not hard to obtain the remaining estimates. We omit the details.
\end{proof}

\begin{lemma}\label{lem:modifed-HY-ineqaulity}
	Under the same assumptions of  Lemma \ref{lem:u-modifiedG}, one has
	\begin{equation}\label{eq:source-inequality}
		\int_MQ_g\ud\mu_g+
		\frac12\int_M|E_{u^2g}|_{u^2g}^2\ud\mu_{u^2g}
		\leq		\frac1{24}\int_M\left(\frac f u\right)^2\ud\mu_g.
	\end{equation}
For any constant $\varepsilon>0$, there holds
	\begin{equation}\label{eq:weighted-scalar}
		\int_MR_g^2\frac{\varepsilon u^2}{(u+\varepsilon)^3}\ud\mu_g
		\leq\frac94\int_M\frac{f^2}{(u+\varepsilon)^2}\ud\mu_g.
	\end{equation}
\end{lemma}

\begin{proof}
	Consider the vector field $\mathcal{K}_\eps(u)$ coming from \eqref{eq:K-eps}.  Using Lemma \ref{lem:u-modifiedG}, one has $\mathcal{K}_\eps(u)\in L^1(M,g)$. Set $\theta=\frac{\varepsilon}{u+\varepsilon}$.
	Multiplying $\mathcal{K}_\eps(u)$ by the  cut-off functions as in \eqref{cut-eta_j} and then taking limit, Lemma \ref{lem:divK} yields
	\begin{align*}
&\int_M\left(1-\theta\right)^2\bigg[Q_g
		+\frac12\big|E_{u^2g}-4\theta\ud \log u\otimes\ud \log u+\theta |\nabla \log u|^2g\big|_g^2\\
		&+\frac{\theta R_g^2}{18}+\frac{\theta R_g}{3}(1+3\theta)|\nabla\log u|^2\bigg]\dmu\\
		=&\int_M\left(1-\theta\right)^2\left(\frac{\theta  R_gL_gu}{18u}+\frac{(L_gu)^2}{24u^2}\right)\dmu\\
		=&\int_M\left(1-\theta\right)^2\left(\frac{\theta  R_gf}{18u}+\frac{f^2}{24u^2}\right)\dmu.
	\end{align*}
Since $f$ has compact support, letting $\eps\to 0$ and using Fatou's lemma, we prove \eqref{eq:source-inequality}. With the help of Young's inequality, one has 
	\begin{align*}
		&\int_M\left(1-\theta\right)^2\left(\frac{\theta  R_gf}{18u}+\frac{f^2}{24u^2}\right)\dmu\\
		\leq &\int_M\left(1-\theta\right)^2\left(\frac{\theta  }{18}\left(\frac{R_g^2}{3}+\frac{3f^2}{4u^2}\right)+\frac{f^2}{24u^2}\right)\dmu
	\end{align*}
Combining with the previous identity, we prove \eqref{eq:weighted-scalar}.
\end{proof}

Define the  Yamabe invariant
\begin{equation*}
	Y(M,[g])=\inf_{0\ne\varphi\in C_c^\infty(M)}
	\frac{\displaystyle\int_M(6|\nabla\varphi|_g^2+R_g\varphi^2)\ud\mu_g}
	{\displaystyle\left(\int_M|\varphi|^4\ud\mu_g\right)^{1/2}}.
\end{equation*}
 Since  $R_g>0$, by testing an Aubin-Talenti bubble  in a small neighborhood of    some point, a direct computation yields
 \begin{equation}\label{Yamabe-range}
 	0\leq Y(M,[g])\leq Y(\ms^4,[g_{\ms^4}])=8\sqrt6\pi.
 \end{equation}
\begin{proposition}\label{prop:yamabepositive}
 There holds
\begin{equation*}
	Y(M,[g])>0.
\end{equation*}
\end{proposition}

\begin{proof}
	We argue by contradiction and suppose that  $Y(M,[g])=0$. 
 Choose a smooth connected exhaustion $\Omega_j$ and set  $Y_j=Y(\Omega_j,[g])$ which tends to zero as $j\to\infty$. For sufficiently large $j$, using \cite[Theorem 3.1]{ACM},  there exists  a positive solution $v_j$ to 
		\begin{equation}\label{eq:dirichlet-yamabe}
		L_gv_j=Y_jv_j^3,\qquad v_j|_{\partial\Omega_j}=0,
		\qquad \int_{\Omega_j} v_j^4\ud\mu_g=1.
	\end{equation}
	Let $w_j=\sqrt{Y_j}v_j$. Then, $$L_gw_j=w_j^3, \quad \|w_j\|_{L^4(\Omega_j)}^2=Y_j\to0, \quad \mathrm{as}\;\; j\to\infty.$$ 

Fix a ball $B_g(p_0,r_0)$ and choose a nonnegative cutoff function $\eta$ satisfying $\eta\equiv 1$ in $B_g(p_0,r_0)$ and $\eta\equiv 0$ outside $B_g(p_0,2r_0)$.   Multiplying  by $\eta^2w_j^3$ and using   Sobolev inequality, one has
	\begin{equation}\label{eq:yamabe-bootstrap}
		\|\eta w_j^2\|_{L^4}^2
		\leq C\|w_j\|_{L^4(\mathrm{supp}(\eta))}^2\|\eta w_j^2\|_{L^4}^2
		+C\int_{\mathrm{supp}(\eta)}w_j^4\ud\mu_g.
	\end{equation}
	For large $j$, the first term on right-hand side is absorbed and then $$\|w_j\|_{L^8(B_g(p_0,r_0))}\leq C\sqrt{Y_j}.$$ 
	Applying  Moser iteration to $-\Delta_gw_j\leq w_j^2w_j/6$ gives $$\sup_{B_g(p_0,r_0/2)} w_j\leq C\sqrt{Y_j}\to0, \quad \mathrm{as}\;\; j\to\infty.$$ 
	
		Extend $f_j=Y_jv_j^3$ by zero.  Since $v_j$ vanishes on $\partial\Omega_j$, the derivatives up to second of $v_j^3$  also vanish. Hence,	such extension belongs to  $C^{2,\alpha}$.  Then Lemma \ref{lem:u-modifiedG}
		shows there exists a positive solution $u_j$ to $L_gu_j=f_j$.  The maximum principle yields $u_j\geq v_j$.
Then,  on $B_g(p_0,r_0/2)$, 
	$$\frac{f_j}{u_j}\leq Y_jv_j^2=w_j^2$$ and tend locally uniformly to zero. Now set  $U_j(x)=u_j(x)/u_j(p_0)$ satisfying  $$-6\Delta_gU_j+(R_g-f_j/u_j)U_j=0, \quad U_j(p_0)=1.$$ For sufficiently  large $j$,  the local Harnack inequality therefore gives  a constant $C$  independent of $j$ such that
		\begin{equation}\label{U_j-bound}
			C^{-1}\leq U_j\leq C,\quad  \mathrm{on}\;  B_g(p_0,r_0/4)
		\end{equation}
	Applying \eqref{eq:weighted-scalar} with $u=u_j$ and $\varepsilon=u_j(p_0)$,  we obtain
	\begin{equation*}
		\int_MR_g^2\frac{U_j^2}{(1+U_j)^3}\ud\mu_g
		\leq\frac94\int_M\frac{f_j^2}{(u_j+u_j(p_0))^2}\ud\mu_g
		\leq\frac94Y_j^2\longrightarrow0, \quad \mathrm{as}\; j\to\infty
	\end{equation*}
which contradicts \eqref{U_j-bound} and $R_g>0$.	 Thus,  we prove $Y(M,[g])>0$.
\end{proof}

\begin{lemma}\label{lem:green-bounds}
For each $s>0$, there holds
	\begin{equation}\label{eq:green-weak-L2}
		\Vol_g(\{G_p\geq s\})\leq\frac1{Y(M,[g])^2s^2}.
	\end{equation}
	For every compact neighborhood $K$ of $p$, $$G_p\in L^4(M\setminus K,g).$$ There is a positive constant $C$ such that, if $p\notin B_g(x,r)$,
	\begin{equation}\label{eq:green-moser}
		\sup_{B_g(x,r/2)}G_p
		\leq\frac{C}{\sqrt{Y(M,[g])}\,r}
		\left(\int_{B_g(x,r)}G_p^4\ud\mu_g\right)^{1/4}.
	\end{equation}
In particular, $G_p(x)\to0$ as $d_g(x,p)\to\infty.$
\end{lemma}

\begin{proof}
	For a Dirichlet Green function $G_\Omega$ extended by zero, test its equation by $v=\min\{G_\Omega,s\}$ and obtain $$\int_\Omega(6|\nabla v|^2+R_gG_\Omega v)\ud\mu_g=s.$$ Since $0\leq v\leq G_\Omega$, then  $$Y(M,[g])\|v\|_{L^4(M)}^2\leq s.$$ Since $v=s$ on $\{G_\Omega\geq s\}$ for $s>0$, there holds
	$$\Vol_g(\{G_\Omega\geq s\})\leq\frac1{Y(M,[g])^2s^2}.$$
Exhaustion gives  	 \eqref{eq:green-weak-L2}.  
	
	The exterior bound \eqref{eq:exterior} and the layer-cake formula give
	\begin{align*}
		\int_{M\setminus K}G^4\ud\mu_g\leq  &4 \int^{\sup_{\partial K}G}_0 s^3\Vol_g(\{G>s\})\ud s\\
		\leq &\frac{2}{Y(M,[g])^2}\left(\sup_{\partial  K}G\right)^2 <\infty.
	\end{align*}
	For each constant $a\geq1$, away from $p$,  direct differentiation gives $$L_g(G^a)=(1-a)(R_g+6a|\nabla\log G|^2)G^a\leq0$$  Testing by $\eta^2G^a$ and applying the Yamabe inequality yields
	\begin{equation}\label{eq:green-moser-step}
		Y(M,[g])\|\eta G^a\|_{L^4}^2\leq6\int_M|\nabla\eta|^2G^{2a}\ud\mu_g.
	\end{equation}
Choose $a=2,4,8,\ldots$ and distance cutoff functions between $B_g(x,r)$ and $B_g(x,r/2)$.  The Moser  iteration gives \eqref{eq:green-moser}.
\end{proof}

\begin{lemma}\label{lem:G-smooth-potential}
	Under the same assumptions of  Lemma \ref{lem:u-modifiedG},
  for a  fixed  ball $B_g(p,\frac{r_0}{2})$ containing $\mathrm{supp}f$, one has
  \begin{equation}\label{uequivG}
  	C^{-1}\leq \frac{u}{G_p}\leq C, \quad \mathrm{on}\;\; M\setminus B_g(p,r_0/2).
  \end{equation}
Moreover, for the conformal metric $h=G_p^2g$, one has
			\begin{equation}\label{eq:ratio-integrability}
			\int_{M\setminus  B_g(p,r_0)}\left(\left|\nabla_h^2\left(\log\frac{u}{G_p}\right)\right|_h^2+\left|\nabla_h \left(\log \frac{u}{G_p}\right)\right|_h^4\right)\ud\mu_h<\infty.
		\end{equation}
		\end{lemma}
\begin{proof}
	Following the argument of Lemma \ref{lem:green} and  choosing a sequence of exhaustion $\Omega_j$ of $M$ where $B_g(p,r_0)\Subset\Omega_1$, such that 
	$$L_g u_j=f \quad \mathrm{in}\;\; \Omega_j, \quad u_j=0\quad \mathrm{on}\;\;\partial \Omega_j.$$
Then the maximum principle yields $u_j>0$ in $\Omega_j$. Letting $j\to \infty$, we obtain a minimal positive solution $u$ to $L_g u=f$ in $M$. 	Meanwhile, let $G_{p,\Omega_j}$ denote the related Dirichlet Green function.

Since $u_j$ and $G_{p,\Omega_j}$ converge to $u$ and $G_p$ respectively, their ratios on $\partial B_g(p,r_0/2)$ therefore have positive uniform lower and upper bounds for  sufficiently large $j$. Both functions vanish on $\partial \Omega_j$, so the maximum principle deduces
$$C^{-1}G_{p,\Omega_j}\leq u_j\leq  CG_{p,\Omega_j}$$
which yields  \eqref{uequivG}  by letting $j\to\infty$.

For brevity, set $v=u/G_p$. 
With the help of the conformal covariant property of \eqref{L_g}, one has
$$\Delta_h (v)=0, \quad \mathrm{on}\;\; M\setminus  B_g(p,r_0/2).$$
Using Theorem \ref{thm:CVandHY-for-Q}, one has 
\begin{equation}\label{finite-Ric_h^2}
	\int_{M\setminus  B_g(p,r_0/2)}|\Ric_h|^2_h\ud\mu_h=\int_{M\setminus B_g(p,r_0/2)}|\Ric_h|^2_g\ud\mu_g\leq 32\pi^2.
\end{equation}

All the following derivatives and norms are taken with respect to
$h$. Let $0\leq\eta\leq1$ be smooth and compactly supported in $M\setminus  B_g(p,r_0/2)$ to be chosen later.
Set the notations
$$	A=\int_{M}\eta^4|\nabla^2v|^2\dmuh,\quad 
B=\int_{M}\eta^4|\nabla v|^4\dmuh,\quad 
D_\eta=\int_{M}|\nabla\eta|^4\dmuh.$$

First, $\Delta_hv=0$ gives
$$
\mathrm{div}_h\bigl(v|\nabla v|^2\nabla v\bigr)
=|\nabla v|^4+2v\nabla^2v(\nabla v,\nabla v).
$$
Multiplying by $\eta^4$ and integrating by parts, we obtain
\begin{equation}\label{eq:gradient-identity}
	\begin{split}
		B={}&-2\int_{M}\eta^4v
		\nabla^2v(\nabla v,\nabla v)\dmuh\\
		&-4\int_{M}\eta^3v|\nabla v|^2
		\ip{\nabla\eta}{\nabla v}_h\dmuh.
	\end{split}
\end{equation}
By H\"older's inequality,  \eqref{uequivG} and Young's inequality, one has 
$$
	B\leq CA^{1/2}B^{1/2}
	+CB^{3/4}D_\eta^{1/4}\leq \frac{1}{2}B+CA+CD_\eta
$$
which yields
\begin{equation}\label{eq:B-from-A}
	B\leq CA+CD_\eta.
\end{equation}
Next, the Bochner formula and $\Delta_hv=0$ give
$$
\frac12\Delta_h|\nabla v|^2
=|\nabla^2v|^2+\Ric_h(\nabla v,\nabla v).
$$
Multiplying by $\eta^4$ and integrating by parts yields
$$	A=-\int_{M}\eta^4\Ric_h(\nabla v,\nabla v)\dmuh-2\int_{M}\eta^3
\langle \nabla\eta, \nabla|\nabla v|^2\rangle_h\dmuh.$$
Using \eqref{finite-Ric_h^2} and H\"older's inequality,
$$
	\left|\int_{M}\eta^4\Ric_h(\nabla v,\nabla v)\dmuh\right|
	\leq C B^{1/2}.
$$
Moreover, 
\begin{align*}
	2\left|\int_{M}\eta^3
	\ip{\nabla\eta}{\nabla|\nabla v|^2}_h\dmuh\right|
	&\leq4\int_{M}\eta^3|\nabla\eta|
	|\nabla^2v||\nabla v|\dmuh\\
	&\leq\frac12A
	+8\int_{M}\eta^2|\nabla\eta|^2|\nabla v|^2\dmuh\\
	&\leq\frac12A+8D_\eta^{1/2}B^{1/2}.
\end{align*}
Consequently,
$$
A\leq C B^{1/2}+16D_\eta^{1/2}B^{1/2}.
$$
Combining with \eqref{eq:B-from-A} and Young's inequality, one has 
\begin{equation}\label{eq:local-estimate}
	A+B\leq C(1+D_\eta),
\end{equation}

Set
\begin{equation}\label{def-s_eps}
	s_\varepsilon=\frac{G}{G+\varepsilon},\quad \theta_\eps=\frac{\eps}{G+\eps}, \quad w=\log G.
\end{equation}
Using \eqref{eq:finiteE_G} and \eqref{finite-Ric_h^2},  it follows that
\begin{equation}\label{eq:s-uniform-L4}
	\begin{split}
		\int_{M\setminus\{p\}}|\nabla_gs_\varepsilon|_g^4\dmu
		&=\int_{M\setminus\{p\}}
		s_\varepsilon^4\theta_\varepsilon^4
		|\nabla_gw|_g^4\dmu\\
		&\leq\int_{M\setminus\{p\}}
		s_\varepsilon^2\theta_\varepsilon^2
		|\nabla_gw|_g^4\dmu\leq C.
	\end{split}
\end{equation}
Choose a fixed smooth function $\chi$ such that
$$
0\leq\chi\leq1,\qquad
\chi\equiv0\ \text{in }\overline B_g(p,3r_0/4),\qquad
\chi\equiv1\ \text{in } M\setminus B_g(p,r_0).
$$
Its gradient is supported in a fixed compact subset of $B_g(p,r_0)\setminus B_g(p,3r_0/4)$.
Choosing the distance cutoff functions $\eta_j$ from \eqref{cut-eta_j}, set
$
\eta=\chi s_\varepsilon\eta_j.
$
Then, 
\begin{align*}
	&\int_{M}|\nabla_h\eta|_h^4\dmuh\\
	\leq& C\Bigg(\int_{M\setminus B_g(p,r_0/2)}|\nabla_gs_\varepsilon|_g^4\dmu	+\int_{M\setminus B_g(p,r_0/2)} s_\varepsilon^4|\nabla_g\eta_j|_g^4\dmu
	\Bigg)+C.
\end{align*}
The first term  is uniformly bounded by \eqref{eq:s-uniform-L4}. For the second term,
use $s_\varepsilon\leq G/\varepsilon$ and
$|\nabla_g\eta_j|_g\leq C/j$ to obtain
\begin{equation}\label{eq:outer-cutoff-error}
	\int_{M\setminus B_g(p,r_0/2)} s_\varepsilon^4|\nabla_g\eta_j|_g^4\dmu
	\leq\frac{C}{j^4\varepsilon^4}\int_{M\setminus B_g(p,r_0/2)} G^4\dmu
\end{equation}
which is finite by using Lemma \ref{lem:green-bounds}. Thus, for
every fixed $\varepsilon>0$, the right-hand side of
\eqref{eq:outer-cutoff-error} tends to zero as $j\to\infty$.

Applying \eqref{eq:local-estimate}  gives
\begin{align*}
	&\int_M\chi^4s_\varepsilon^4\eta_j^4
	\left(|\nabla_h^2v|_h^2+|\nabla_hv|_h^4\right)\dmuh\\
	&\hspace{15mm}\leq
	C+\frac{C}{j^4\varepsilon^4}\int_{M\setminus B_g(p,r_0/2)} G^4\dmu,
\end{align*}
where $C$ is independent of $j$ and $\varepsilon$.
First fix $\varepsilon>0$ and let $j\to\infty$.
Fatou's lemma gives
$$
\int_M\chi^4s_\varepsilon^4
\left(|\nabla_h^2v|_h^2+|\nabla_hv|_h^4\right)\dmuh\leq C.
$$
Since $G>0$,
$s_\varepsilon\to1$ as $\eps\to 0$ on $M\setminus B_g(p,r_0/2)$. A second application of
Fatou's lemma yields
$$
\int_{M}\chi^4
\left(|\nabla_h^2v|_h^2+|\nabla_hv|_h^4\right)\dmuh\leq C.
$$
 In $M\setminus B_g(p,r_0/2)$, the estimate \eqref{uequivG}  yields
$$
|\nabla_h\log v|_h^4\le C|\nabla_hv|_h^4,
\qquad
|\nabla_h^2\log v|_h^2
\le
C|\nabla_h^2v|_h^2
+C|\nabla_hv|_h^4.
$$
Thus,  we prove \eqref{eq:ratio-integrability}.
\end{proof}

\begin{proposition}\label{prop:total-q}
	Under the same assumptions of  Lemma \ref{lem:u-modifiedG}, one has 
	$$
	\Ric_{u^2g}\in L^2(M,u^2g),
	\qquad
	Q_{u^2g}\in L^1(M,u^2g).
	$$
	Moreover, one has 
	\begin{equation}\label{eq:total-q}
		\int_M Q_{u^2g}\,\ud\mu_{u^2g}
		=
		16\pi^2-\frac12
		\int_{M\setminus\{p\}}|\Ric_{G_p^2g}|_g^2\,\ud\mu_g.
	\end{equation}
\end{proposition}

\begin{proof}
For brevity, set 
	$$
	v=\frac{u}{G}, \quad h=G^2g.
	$$
	The conformal covariance of the conformal Laplacian gives
	$$
	R_h=G^{-3}L_gG=0
	\quad\text{in }M\setminus\{p\},
	\qquad
	R_{u^2g}=u^{-3}L_gu=\frac{f}{u^3}.
	$$
	In particular, $R_{u^2g}\in C_c^{2,\alpha}(M)$. Meanwhile,
	$$
	-6\Delta_hv=L_hv=G^{-3}f.
	$$
	It follows that
	\begin{equation*}
		\Delta_h\log v	=
		-|\nabla_h\log v|_h^2-\frac{f}{6uG^2}
		\qquad\text{in }M\setminus\{p\}.
	\end{equation*}
Since $\mathrm{supp} f\subset B_g(p,r_0/2)$, one has
\begin{equation}\label{eq:log-equation}
	\Delta_h\log v=-|\nabla_h\log v|_h^2
	\qquad\text{in }M\setminus B_g(p,r_0).
\end{equation}
Since 
	\begin{equation*}
		\Ric_{u^2g}
		=
		\Ric_h-2\nabla_h^2\log v+2\,\ud \log v\otimes\ud \log v
		-\bigl(\Delta_h\log v+2|\nabla_h\log v|_h^2\bigr)h.
		\end{equation*}
 in $M\setminus B_g(p,r_0)$, one has 
	$$
	|\Ric_{u^2g}|_h^2
	\le C\left(
	|\Ric_h|_h^2+|\nabla_h^2\log v|_h^2+|\nabla_h\log v|_h^4
	\right).
	$$
Combining with \eqref{finite-Ric_h^2} and \eqref{eq:ratio-integrability}, we obtain
\begin{align*}
	\int_{M\setminus B_g(p,r_0)}
	|\Ric_{u^2g}|_{u^2g}^2\,\ud\mu_{u^2g}
	=	\int_{M\setminus B_g(p,r_0)}
	|\Ric_{u^2 g}|_{h}^2\,\ud\mu_{h}\leq C.
\end{align*}
	Since $u$ is positive and regular near $p$, then
	\begin{equation}\label{eq:Ric-tilde-L2}
		\int_M|\Ric_{u^2g}|_{u^2 g}^2
		\ud\mu_{u^2 g}<\infty.
	\end{equation}
Recalling
	$R_{u^2g}=f/u^3\in C_c^{2,\alpha}(M)$, \eqref{eq:Ric-tilde-L2} deduces
	\begin{equation}\label{eq:Q-tilde-L1}
		Q_{u^2g}\in L^1(M,u^2g).
	\end{equation}
	
Since $R_{h}(x)=0$ for  $x\not=p$, one has $Q_h=-\frac{1}{2}|\Ric_h|^2_h.$ Then  \eqref{Paneitz-operator} gives
	\begin{equation}\label{eq:relative-q}
		P_h(\log v)\ud\mu_h
		=
		Q_{u^2g}\,\ud\mu_{u^2g}
		+\frac12|\Ric_h|_h^2\,\ud\mu_h, \quad \text{in} \;\;M\setminus\{p\}.
	\end{equation}
	With the help of  \eqref{finite-Ric_h^2} and \eqref{eq:Q-tilde-L1}, one has 
	\begin{equation}\label{eq:PF-L1}
		P_h(\log v)\in L^1(M\setminus\{p\},h).
	\end{equation}
	
Consider  the  Paneitz vector field  from \eqref{eq:Paneitz-vector} 
	$$
	\mathcal{J}_h(\log v)=
	\nabla_h\Delta_h\log v+2\Ric_h^\sharp(\nabla_h\log v)
	\qquad\text{in }M\setminus\{p\}.
	$$
	In $M\setminus B_g(p,r_0)$, using \eqref{eq:log-equation}, there holds
	\begin{align*}
		\mathcal J_h(\log v)
		&=-\nabla_h|\nabla_h\log v|_h^2+2\Ric_h^\sharp(\nabla_h\log v)\\
		&=-2(\nabla_h^2\log v)^\sharp(\nabla_h\log v)
		+2\Ric_h^\sharp(\nabla_h\log v).
	\end{align*}
Combining with H\"older's inequality and \eqref{eq:ratio-integrability}, we have 
$$
\|\mathcal J_h(\log v)\|_{L^{\frac{4}{3}}(M\setminus B_g(p,r_0),h)}\leq C
$$
 Since $\mathcal{J}_h(\log v)$
	is regular on every compact set away  from $p$, it follows that,
	for each sufficiently small fixed constant $\delta>0$,
	\begin{equation}\label{eq:J-L43}
		\mathcal J_h(\log v)\in L^{4/3}(M\setminus B_g(p,\delta),h).
	\end{equation}
Set $s_\eps,\theta_\eps, w$ same as in \eqref{def-s_eps}.
	Then
	$$
	s_\eps+\theta_\eps=1,
	\qquad
	\ud s_\eps=s_\eps\theta_\eps\,\ud w,
	\qquad
	\left|\ud w\otimes\ud w-\frac14|\ud w|_g^2g\right|_g^2=\frac34|\ud w|_g^4.
	$$
	Consequently, using \eqref{eq:finiteE_G}, \eqref{finite-Ric_h^2} and $0\leq s_\eps\leq 1$,   one has 
	\begin{align*}
		12\int_{M\setminus\{p\}}
		s_\eps^2\theta_\eps^2|\ud w|_g^4\,\ud\mu_g
		&=
		\int_{M\setminus\{p\}}s_\eps^2\left|4\theta_\eps \ud w\otimes\ud w-\theta_\eps|\ud w|_g^2g\right|_g^2\,\ud\mu_g\\
		&\le
		2\int_{M\setminus\{p\}}
		s_\eps^2\left|\Ric_h-4\theta_\eps \ud w\otimes\ud w+\theta_\eps|\ud w|_g^2g\right|_g^2\,\ud\mu_g\\
		&\quad+
		2\int_{M\setminus\{p\}}
		s_\eps^2|\Ric_h|_g^2\,\ud\mu_g\leq C.
	\end{align*}
	Then, using $s_\eps\leq 1, \theta_\eps\le1$, we obtain
	\begin{equation}\label{eq:cutoff-L4}
		\int_{M\setminus\{p\}}|\nabla_h s_\eps|_g^4\,\ud\mu_g=
	\int_{M\setminus\{p\}}
	s_\eps^4\theta_\eps^4|\ud w|_g^4\,\ud\mu_g\leq C.
	\end{equation}
	Choose the cutoff functions $\eta_j$ as in  \eqref{cut-eta_j}.
	For $j$ sufficiently large,  multiplying $P_h(\log v)$ by $s_\eps\eta_j$ and  integrating  by parts over $M\setminus B_g(p,\delta)$, one has 
	\begin{equation}\label{eq:cutoff-identity}
		\begin{aligned}
			&\int_{M\setminus B_g(p,\delta)}s_\eps\eta_jP_h(\log v)\,\ud\mu_h\\
			=&
			\int_{\partial B_g(p,\delta)}
			s_\eps\langle\mathcal J_h(\log v),\nu_h\rangle_h\,\ud\sigma_h\\
			&-\int_{M\setminus B_g(p,\delta)}
			\eta_j\langle\mathcal J_h(\log v),\nabla_hs_\eps\rangle_h\,\ud\mu_h\\
			&-\int_{M\setminus B_g(p,\delta)}
			s_\eps\langle\mathcal J_h(\log v),\nabla_h\eta_j\rangle_h\,\ud\mu_h,
		\end{aligned}
	\end{equation}
	where $\nu_h$ is the outward unit normal of $M\setminus B_g(p,\delta)$;
	thus it points into the  ball.
	
	We first let $j\to\infty$ with $\eps$ and $\delta$ fixed.
Using $s_\eps\le G/\eps$ and Lemma \ref{lem:green-bounds}, one has 
	\begin{align*}
		\int_{M\setminus B_g(p,\delta)}
		s_\eps^4|\nabla_h\eta_j|_h^4\,\ud\mu_h
		&=
		\int_{M\setminus B_g(p,\delta)}
		s_\eps^4|\nabla_g\eta_j|_g^4\,\ud\mu_g\\
		&\le
		\frac{C}{\eps^4j^4}
		\int_{M\setminus B_g(p,r_0)}G^4\,\ud\mu_g
		\to 0,  \quad \mathrm{as}\;\; j\to\infty.
	\end{align*}
	It follows from H\"older's inequality and \eqref{eq:J-L43} that
	\begin{align*}
			&\left|
		\int_{M\setminus B_g(p,\delta)}
		s_\eps\langle\mathcal J_h(\log v),\nabla_h\eta_j\rangle_h\,\ud\mu_h
		\right|\\
		\leq &C\|s_\eps\nabla_h\eta_j\|_{L^4(M\setminus B_g(p,\delta),h)}
		\longrightarrow0, \quad \mathrm{as}\;\; j\to\infty.
	\end{align*}
Similar, using \eqref{eq:J-L43} and \eqref{eq:cutoff-L4}, one has
$$\left|\int_{M\setminus B_g(p,\delta)}
\eta_j\langle\mathcal J_h(\log v),\nabla_hs_\eps\rangle_h\,\ud\mu_h\right|\leq C.$$
 Combining with \eqref{eq:PF-L1} and applying  the dominated convergence theorem in  \eqref{eq:cutoff-identity} for $j\to\infty$, one has 
	\begin{equation}\label{eq:epsilon-identity}
		\begin{aligned}
			\int_{M\setminus B_g(p,\delta)}s_\eps P_h(\log v)\,\ud\mu_h
			=&
			\int_{\partial B_g(p,\delta)}
			s_\eps\langle\mathcal J_h(\log v),\nu_h\rangle_h\,\ud\sigma_h\\
			&-\int_{M\setminus B_g(p,\delta)}
			\langle\mathcal J_h(\log v),\nabla_hs_\eps\rangle_h\,\ud\mu_h.
		\end{aligned}
	\end{equation}
	
	Next, we claim 
	\begin{equation}\label{eq:epsilon-error}
		\int_{M\setminus B_g(p,\delta)}
		\langle\mathcal J_h(\log v),\nabla_hs_\eps\rangle_h\,\ud\mu_h
		\longrightarrow0
		\qquad\text{as }\eps\to 0.
	\end{equation}
	Given $\eps_1>0$, choose $R$ sufficiently large that
	$$
	\|\mathcal J_h(\log v)\|_{L^{4/3}(M\setminus B_g(p,R),h)}<\eps_1.
	$$
	Then,  using H\"older's inequality, \eqref{eq:cutoff-L4} implies
	$$
	\left|
	\int_{M\setminus B_g(p,R)}
	\langle\mathcal J_h(\log v),\nabla_hs_\eps\rangle_h\,\ud\mu_h
	\right|
	\le C\eps_1.
	$$
On the compact annulus
	$\overline{B_g(p,R)\setminus B_g(p,\delta)}$,
	$\nabla_hs_\eps\to0$ uniformly, so the corresponding integral
	tends to zero. Taking  $\eps\to 0$ and
	then letting $\eps_1\to 0$ proves the claim \eqref{eq:epsilon-error}.
	
Letting  $\eps\to0$ in
	\eqref{eq:epsilon-identity} and using \eqref{eq:PF-L1}, the dominated
	convergence yields 
	\begin{equation}\label{eq:inner-flux}
		\int_{M\setminus B_g(p,\delta)}P_h(\log v)\,\ud\mu_h
		=
		\int_{\partial B_g(p,\delta)}
		\langle\mathcal J_h(\log v),\nu_h\rangle_h\,\ud\sigma_h
	\end{equation}
where $\nu_h$ is the inner normal vector of $\partial B_g(p,\delta)$ with respect to $h$.
	
	Using Lemma \ref{lem:Green-behavior},	after a delicate but standard computation, one has 
	$$
	\lim_{\delta\to0}
	\int_{\partial B_g(p,\delta)}
	\langle\mathcal J_h(\log v),\nu_h\rangle_h\,\ud\sigma_h
	=16\pi^2.
	$$
	
	Finally, letting $\delta\to0$ in \eqref{eq:inner-flux} and using \eqref{eq:PF-L1}, there holds
	$$
	\int_{M\setminus\{p\}}P_h(\log v)\,\ud\mu_h=16\pi^2.
	$$
	Integrating \eqref{eq:relative-q} and using
	$
	|\Ric_h|_h^2\,\ud\mu_h
	=
	|\Ric_{G_p^2g}|_g^2\,\ud\mu_g
	$, we  prove \eqref{eq:total-q}.
\end{proof}

\section{Improved Hang-Yang-type inequality}\label{sec:impoved-HY}

Throughout this section, we let $(M,g)$ be an open four-manifolds with 
 $Q_g\geq 0$ and $R_g>0$.
 
Using Theorem  \ref{thm:CVandHY-for-Q}, we define the following nonnegative quantity
\begin{equation}\label{def-B_g}
	\mathcal{B}_g=16\pi^2-\frac12\int_{M\setminus\{p\}}
	|\Ric_{G_p^2g}|_{g}^2\ud\mu_{g}
	-\int_M Q_g\ud\mu_g
\end{equation}
which is independent  of $p$ due to  Proposition \ref{prop:total-q}.
\begin{theorem}\label{thm:yamabe-total-defect}
	 There holds
	\begin{equation}\label{eq:source-yamabe}
		\int_MQ_g\ud\mu_g+\mathcal{B}_g\leq\frac{Y(M,[g])^2}{24}.
	\end{equation}
\end{theorem}

\begin{proof}
	If $Y(M,[g])=Y(\mathbb S^4)$, \eqref{eq:final-g} shows that  \eqref{eq:source-yamabe} holds automatically.

If $Y(M,[g])<Y(\mathbb S^4)$. Choose a smooth connected exhaustion $\Omega_j$ with $Y_j=Y(\Omega_j,[g])\to Y(M,[g])$. Same as before, for all sufficiently large $j$, there exists a  positive minimizer $v_j$, normalized by $\int_{\Omega_j}v_j^4\ud\mu_g=1$, satisfying $L_gv_j=Y_jv_j^3$ and $v_j=0$ on $\partial\Omega_j$ (See \cite[Theorem 3.1]{ACM}).
	
	Extend $f_j=Y_jv_j^3$ by zero and let $u_j$ be the positive minimal solution to  $L_gu_j=f_j$. The maximum principle gives $u_j\geq v_j$ on $\Omega_j$, and therefore
	\begin{equation}\label{R_uj^2}
	\int_MR_{u_j^2g}^2\ud\mu_{u_j^2g}
		=\int_{\Omega_j}\frac{f_j^2}{u_j^2}\ud\mu_g
		\leq Y_j^2\int_{\Omega_j}v_j^4\ud\mu_g=Y_j^2.
	\end{equation}
Since $R_{u_j^2g}$ has a compact support, applying  Proposition \ref{prop:total-q} and \eqref{R_uj^2} ,  one has
	\begin{align*}
		\int_MQ_g\dmu+\mathcal{B}_g=&16\pi^2-\frac{1}{2}\int_{M\setminus\{p\}}|\Ric_{G^2g}|^2_g\dmu\\
		=&\int_MQ_{u_j^2g}\ud\mu_{u_j^2g}\\
	= &-\frac{1}{2}\int_M|E_{u_j^2g}|_{u_j^2g}^2\ud\mu_{u_j^2g}+\frac{1}{24}\int_MR_{u_j^2g}\ud\mu_{u_j^2g}\leq\frac{1}{24}Y_j^2.
	\end{align*}
Letting $j\to\infty$, we prove \eqref{eq:source-yamabe}.
\end{proof}

During the remaining of this section, we set the following notations
$$w(x)=-\log G_p(x), \quad  \Omega_t=\{p\}\cup\{w(x)\leq t\}, \quad \Sigma_t=\{x:w=t\}.$$
Then, at  the  regular value of $w$, we set 
\begin{equation*}
	\begin{aligned}
		a(t)&=\int_{\Sigma_t} |\nabla w|\ud\sigma_g,
		&m(t)&=\int_{\Sigma_t}\frac1{|\nabla w|}\ud\sigma_g,\\
		D(t)&=\int_{\Sigma_t}\left(|\nabla w|^2+\frac{R_g}{6}\right)|\nabla w|\ud\sigma_g,
		&F(t)&=\int_{\Sigma_t}\frac{(|\nabla w|^2+R_g/6)^2}{|\nabla w|}\ud\sigma_g.
	\end{aligned}
\end{equation*}

Since $\Delta_gG_p=R_gG_p/6>0$ on $M\setminus\{p\}$ and the Sard's theorem, the critical set of $G_p$ has zero measure so  is $w$. Thus, the critical values have no essential  impact on the following lemmas.

\begin{lemma}\label{lem:level-identities}
	For almost every regular value, the level functions satisfy
	\begin{equation}\label{a'}
		a'(t)=a(t)-\int_{\Sigma_t}\frac{R_g}{6|\nabla w|}\ud\sigma_g\leq a(t);
	\end{equation}
\begin{equation}\label{2a-a'}
	(2a-a')^2\leq  mF; 
\end{equation}
	\begin{equation}\label{eq:level-a-pole}
		e^{-t}a(t)\to\frac{1}{6}\quad \mathrm{as}\;\;  t\to-\infty; 
		\qquad 0<a(t)\leq\frac16e^t.
	\end{equation}
	For each $\xi\in C^1_c(\mr)$, one has 
	\begin{equation}\label{eq:level-yamabe}
		Y(M,[g])\left(\int_\mr m(t)\xi(t)^4\ud t\right)^{1/2}
		\leq6\int_\mr a(t)(\xi'(t)+\xi(t))^2\ud t.
	\end{equation}
For each  fixed $T_0$ and $\varepsilon>0$, there is a constant $C_{\varepsilon,T_0}$ such that,  for $T\geq T_0$, 
	\begin{equation}\label{eq:F-average}
		\int_{T_0}^T F(t)\ud t\leq(\mathcal{B}_g+\varepsilon)(T-T_0)+C_{\varepsilon,T_0}.
	\end{equation}
\end{lemma}

\begin{proof}
Using  divergence theorem,  one has
$$
	a(t_2)-a(t_1)=\int_{\{x:t_1<t< t_2\}}\Delta_gw=\int_{t_1}^{t_2}\int_{\Sigma_t}\frac{\Delta_g w}{|\nabla w|}\ud\sigma_g\ud t
$$
Since $\Delta_gw=|\nabla w|^2-R_g/6$, we prove \eqref{a'}.
Consequently, one has 
$$2a-a'=\int_{\Sigma_t}\frac{|\nabla w|^2+R_g/6}{|\nabla w|}\ud\sigma_g.$$
Then,  H\"older's equality yields \eqref{2a-a'}.

Using Lemma \ref{lem:Green-behavior},   one has   
$$e^{-t}a(t)\to \frac{1}{6}, \quad \mathrm{as}\;\; t\to-\infty.$$

Integrating $L_gG=\delta_p$ over $\Omega_t$ gives
	\begin{equation}\label{eq:e^-ta}
		6e^{-t}a(t)+\int_{\Omega_t}R_gG\ud\mu_g=1,
	\end{equation}
	which proves 
	$$a(t)\leq \frac{e^t}{6}$$
and shows that 
\begin{equation}\label{montone-a}
	e^{-s}a(s)\geq e^{-t}a(t), \quad s\leq t.
\end{equation}
 Thus $$a(t)>0$$ for all $t\in \mr$ otherwise $a(t_0)=0$ yields $a(t)=0$ for all $t\geq t_0$ which contradicts $\{x\in M :|\nabla w(x)|=0\}$ has zero measure. Thus, we  prove \eqref{eq:level-a-pole}.
	
Using the coarea formula,  \eqref{a'} and integrating by parts, one has 
\begin{align*}
	&\int_M\bigl(6|\nabla(\xi(w(x))|^2+R_g\xi(w(x))^2\bigr)\ud\mu_g\\
	&=6\int_\mr\bigl(a\xi'^2+(a-a')\xi^2\bigr)\ud t\\
	&=6\int_\mr a(\xi'+\xi)^2\ud t
\end{align*}
which  proves \eqref{eq:level-yamabe}.

	Let $h=G^2g$ and outward normal $\nu=\frac{\nabla w}{|\nabla w|}$.
The conformal Ricci formula \eqref{conformal-Ricci} yields
\begin{equation}\label{Ric-nablaw-nu}
	\Ric_h(\nabla w,\nu)=\Ric_g(\nabla w,\nu)+\partial_\nu(|\nabla w|^2)+(|\nabla w|^2-\frac{R_g}{6})|\nabla w|.
\end{equation}
Using the Paneitz vector field \eqref{eq:Paneitz-vector},   one has 
	\begin{equation*}\label{eq:level-current}
		-\mathcal{J}_g(w)\cdot\nu
		=\partial_\nu\left(|\nabla w|^2+\frac{R_g}{6}\right)
		+2\left(|\nabla w|^2+\frac{R_g}{6}\right)|\nabla w|
		-2\Ric_h(\nabla w,\nu).
	\end{equation*}

	Differentiating $D(t)$ by coarea formula  gives
	$$D'(t)+F(t)=\int_{\Sigma_t}\left(\partial_\nu(|\nabla w|^2+\frac{R_g}{6})+2(|\nabla w|^2+\frac{R_g}{6})|\nabla w|\right)\ud\sigma_g.$$
Combining with \eqref{Ric-nablaw-nu}, one has
\begin{equation}\label{D'+F}
	\begin{split}
		&D'(t)+F(t)-2\int_{\Sigma_t}\Ric_h(\nabla w,\nu)\ud\sigma_g\\
		=&\int_{\Sigma_t}\left(\partial_\nu(-|\nabla w|^2+\frac{R_g}{6})+\frac{2}{3}R_g|\nabla w|-2\Ric_g(\nabla w,\nu)\right)\ud\sigma_g\\
		=&-\int_{\Sigma_t}\langle \mathcal{J}_g(w),\nu\rangle\ud\sigma_g\\
		=&16\pi^2-\int_{\Omega_t}(Q_g+\frac{1}{2}|\Ric_h|^2_g)\dmu.
	\end{split}
\end{equation}
Using \eqref{def-B_g} and $Q_g\geq 0$,  for sufficiently large $T_\eps$, one has for any $T\geq T_\eps$, one has
$$16\pi^2-\int_{\Omega_T}(Q_g+\frac{1}{2}|\Ric_h|^2_g)\dmu \leq \mathcal{B}_g+\frac{\eps}{2}.$$

Then, using \eqref{finite-Ric_h^2}, H\"older's inequality and Young's inequality, we obtain
\begin{align*}
	&2\int_{T_0}^{T}\int_{\Sigma_t}
	\frac{\Ric_h(\nabla w,\nabla w)}{|\nabla w|}\ud\sigma_g\ud t\\
	=&2\int_{\{T_0<w<T\}}\Ric_h(\nabla w,\nabla w)\dmu\\
\leq&2\left(\int_{\{w>T_0\}}|\Ric_h|_g^2\ud\mu_g\right)^{1/2}
	\left(\int_{\{T_0<w<T\}}|\nabla w|^4_g\ud\mu_g\right)^{1/2}\\
\leq & 2\left(\int_{\{w>T_0\}}|\Ric_h|_g^2\ud\mu_g\right)^{1/2}
\left(\int_{\{T_0<w<T\}}(|\nabla w|_g^2+R_g/6)^2\ud\mu_g\right)^{1/2}\\
=&2\left(\int_{\{w>T_0\}}|\Ric_h|_g^2\ud\mu_g\right)^{1/2}\left(\int^{T}_{T_0}F(t)\ud t\right)^{1/2}\leq \frac{\eps}{2(\mathcal{B}_g+\eps)}\int_{T_0}^TF(t)\ud t+C(\eps)
\end{align*}
Integrating \eqref{D'+F} from $T_0$ to $T$ and using  the fact $D(T)\geq0$,  one has
$$D(T)-D(T_0)+\frac{\mathcal{B}_g+\eps/2}{\mathcal{B}_g+\eps}\int_{T_0}^TF(t)\ud t\leq (\mathcal{B}_g+\eps/2)(T-T_\eps)+C(\eps)$$
which proves \eqref{eq:F-average}.
\end{proof}

\begin{theorem}\label{thm:defect-bound}
	Let $(M,g)$ be an open four-manifold with $Q_g\geq 0$ and $R_g>0$. Then
	\begin{equation}\label{eq:sharp-defect}
		\mathcal{B}_g\geq\frac{Y(M,[g])^2}{48}>0.
	\end{equation}
\end{theorem}
\begin{proof}
	First, using Proposition \ref{prop:yamabepositive}, one has $Y(M,[g])>0$.
Define two new functions
\begin{equation}\label{eq:capacity-coordinate}
	\tau(t)=\frac12\log\left(\int_{-\infty}^t\frac{e^{2s}}{a(s)}\ud s\right),
	\qquad b(t)=a(t)e^{2\tau(t)-2t}.
\end{equation}
Using \eqref{eq:level-a-pole}, one has
\begin{equation*}
	e^{2\tau(t)}\geq 6\int^t_{-\infty}e^s\ud s=6e^{t}
\end{equation*}
which yields
$$\tau(t)\geq \frac{t}{2}+\frac{\log 6}{2}.$$
Combining with $e^{-t}a(t)\to \frac{1}{6}$ as $t\to-\infty$, one has  $\tau(t)- (\frac{t}{2}+\frac{\log 6}{2})\to 0$ as $t\to-\infty$. Thus, $\tau(t)$ is an increasing $C^1$ bijection from $\mr$ onto $\mr$. Using \eqref{montone-a}, we have
\begin{equation}\label{b-range}
0<	b(t)=\frac{a(t)}{e^{2t}}\int^t_{-\infty}\frac{e^{2s}}{a(s)}\ud s\leq \int^t_{-\infty}e^{s-t}\ud s=1.
\end{equation}
 Integration and differentiation yield
\begin{equation}\label{eq:capacity-properties}
	\tau'=\frac1{2b},
	\qquad b'=1-\left(2-\frac{a'}a\right)b.
\end{equation}

Let $t=t(\tau)$ be the inverse and set
\begin{equation}\label{eq:capacity-density}
	n(\tau)=\frac{2b(t(\tau))^3m(t(\tau))}{a(t(\tau))^2}.
\end{equation}
Then the volume measure satisfies 
\begin{equation}\label{eq:transformed-volume}
	m(t)\ud t=n(\tau)e^{4(t(\tau)-\tau)}\ud\tau.
\end{equation}
For each $v\in C_c^1(\mr)$, choose $\xi(t)=e^{\tau(t)-t}v(\tau(t))$. A direct computation yields
	$$\xi'+\xi=e^{\tau-t}\frac{v'+v}{2b}, \quad  \ud t=2b\ud\tau.$$ 
Then, there holds
	\begin{equation*}
		6\int_\mr a(\xi'+\xi)^2\ud t
		=3\int_\mr(v'+v)^2\ud\tau
		=3\int_\mr(v'^2+v^2)\ud\tau,
		\quad \int_\mr m\xi^4\ud t=\int_\mr nv^4\ud\tau.
	\end{equation*}
 Using \eqref{eq:level-yamabe}, we obtain
\begin{equation}\label{eq:transformed-yamabe}
	Y(M,[g])\left(\int_\mr n(\tau)v(\tau)^4\ud\tau\right)^{1/2}
	\leq3\int_\mr(v'(\tau)^2+v(\tau)^2)\ud\tau.
\end{equation}

	A direct computation yields
	\begin{equation}\label{intsech}
		\int_\mr\left((\mathrm{sech}'(\tau))^2+\mathrm{sech}^2(\tau)\right)\ud \tau=\frac{8}{3}, \quad \int_\mr \mathrm{sech}^4(\tau)\ud\tau =\frac{4}{3}.
	\end{equation}
	Applying  \eqref{eq:transformed-yamabe} to compact cutoffs of $\mathrm{sech}(\tau-s)$ and 
passing to the limit by Fatou's lemma,  we obtain
	\begin{equation}\label{eq:sech-test}
		Y(M,[g])^2\int_\mr n(\tau)\mathrm{sech}^4(\tau-s)\ud\tau\leq64, 
		\qquad \forall s\in\mr.
	\end{equation}
Immediately,  for all $A\in \mr$, choosing $s=A+1/2$, one has
\begin{equation}\label{unit-int-n}
	\int_A^{A+1} n(\tau)\ud\tau\leq\frac{64\cosh^4(1/2)}{Y(M,[g])^2}.
\end{equation}
Integrating  \eqref{eq:sech-test} over $s\in[\tau(T_0)-K,\tau(T)+K]$ and using Fubini theorem, we obtain
\begin{align*}
	\frac{64(\tau(T)-\tau(T_0)+2K)}{Y(M,[g])^2}\geq &\int^{\tau(T)+K}_{\tau(T_0)-K}\int_{\tau(T_0)}^{\tau(T)} n(\tau)\mathrm{sech}^4(\tau-s)\ud\tau\ud s\\
	\geq &\int^{\tau(T)}_{\tau(T_0)}n(\tau)\ud \tau\int^{K}_{-K}\mathrm{sech}^4(s)\ud s
\end{align*}
Using \eqref{intsech} and  choosing $K$  large enough,  we prove  
		\begin{equation}\label{eq:density-upper}
		\int_{\tau(T_0)}^{\tau(T)}n(\tau)\ud\tau
		\leq\left(\frac{48}{Y(M,[g])^2}+\varepsilon\right)(\tau(T)-\tau(T_0))+C_\varepsilon,
		\qquad \forall T\geq T_0.
	\end{equation}

	Combining \eqref{2a-a'} and \eqref{eq:capacity-properties}, one has 
	 $$F\cdot \frac{b^2m}{a^2}\geq(1-b')^2.$$ Together with H\"older's inequality, there holds
	 $$	\int_{\tau(T_0)}^{\tau(T)}n(\tau)\ud\tau=\int_{T_0}^T\frac{b^2m}{a^2}\ud t\geq\frac{\bigl(T-T_0-b(T)+b(T_0)\bigr)^2}
	 {\int_{T_0}^T F(t)\ud t}.$$
Then, using  $0<b\leq1$ and \eqref{eq:F-average},  one has 
\begin{equation}\label{eq:density-lower}
	\int_{\tau(T_0)}^{\tau(T)}n(\tau)\ud\tau
	\geq\frac{T-T_0}{\mathcal{B}_g+\varepsilon}-C_{\varepsilon,T_0}, \quad 
	\forall T\geq T_0.
\end{equation}

To prove \eqref{eq:sharp-defect}, we argue by contradiction and suppose  that $$48\mathcal{B}_g<Y(M,[g])^2.$$ Choose positive constants  $\varepsilon_1,\varepsilon_2$  small enough such  that 
	$$c:=(\mathcal{B}_g+\varepsilon_1)\cdot \left(\frac{48}{Y(M,[g])^2}+\varepsilon_2\right)<1.$$ Combining \eqref{eq:density-upper} with  \eqref{eq:density-lower}, one has 
	 $$T\leq c\tau(T)+C$$ for all sufficiently large $T$. Hence 
	 $$t(\tau)-\tau\leq-(1-c)\tau+C.$$  Using  coarea formula, \eqref{eq:transformed-volume} and \eqref{unit-int-n}, one has 
\begin{align*}
		\Vol_g(\{w>T\})
			=&\int_{\tau(T)}^\infty n(s)e^{4(t(s)-s)}\ud s\\
			\leq & C\int_{\tau(T)}^\infty n(s)e^{-4(1-c)s}\ud s\\
	\leq &C \sum^\infty_{j=0}e^{-4(1-c)(\tau(T)+j)}\int^{\tau(T)+j+1}_{\tau(T)+j}n(s)\ud s\\
	\leq &C e^{-4(1-c)\tau(T)}\leq Ce^{-\frac{4(1-c)}{c}T}
\end{align*}
which  is equivalent to
\begin{equation}\label{eq:VolG_sbound}
	\Vol_g(\{G_p<s\})\leq Cs^{\frac{4(1-c)}{c}}
\end{equation}
for all  $0<s\leq s_0$  and sufficiently small constant $s_0$.

	Let $x\not= p$ and consider the set $\{p\}\cup\{G_p\geq2G_p(x)\}$ which is compact due to Lemma \ref{lem:green-bounds}. For $r<d_g\bigl(x,\{p\}\cup\{G_p\geq2G_p(x)\}\bigr)$, the ball $B_g(x,r)$ avoids the pole and satisfies $G_p<2G_p(x)$. Applying \eqref{eq:green-moser} at its center and then letting $r\to d_g\bigl(x,\{p\}\cup\{G_p\geq2G_p(x)\}\bigr)$ gives
\begin{equation*}
	d_g\bigl(x,\{p\}\cup\{G_p\geq2G_p(x)\}\bigr)
	\leq\frac{C}{\sqrt{Y(M,[g])}}\Vol_g(\{G_p<2G_p(x)\})^{1/4}.
\end{equation*}
Combining with \eqref{eq:VolG_sbound}, one has
\begin{equation}\label{eq:green-level-distance}
	d_g\bigl(x,\{p\}\cup\{G_p\geq2G_p(x)\}\bigr)
	\leq C(G_p(x))^{\frac{1-c}{c}}.
\end{equation}

Starting from any  point  $y_0$ with $G_p(y_0)<s_0/2$, choose a nearest point $y_1$ on the superlevel  $\{p\}\cup\{G_p\geq2G_p(y_0)\}$. Then, $G_p(y_1)=2G_p(y_0)$ and \eqref{eq:green-level-distance} yields
$$d_g(y_0,y_1)\leq C (G_p(y_0))^{\frac{1-c}{c}}.$$
Then, start at the point  $y_1$ and repeat the same procedure to obtain $y_2, \cdots, y_k$ 
until $G_p(y_{k-1})=2^{k-1}G_p(y_0)<\frac{s_0}{2}$ and $G_p(y_k)=2^kG_p(y_0)\geq s_0/2$.
  In particular, 
$$2^k G_p(y_0)<s_0.$$
Meanwhile, 
$$d_g(y_j,y_{j+1})\leq C(G_p(y_j))^{\frac{1-c}{c}}=CG_p(y_0)^{\frac{1-c}{c}}2^{\frac{1-c}{c}j}.$$
Thus, 
\begin{align*}
d_g(y_0, \{p\}\cup\{G_p\geq s_0/2\})\leq & \sum_{j=0}^{k-1}d_g(y_j, y_{j+1})\\
	\leq &CG_p(y_0)^{\frac{1-c}{c}}\sum^{k-1}_{j=0}2^{\frac{1-c}{c}j}\\
	\leq &CG_p(y_0)^{\frac{1-c}{c}}\frac{2^{\frac{1-c}{c}k}}{2^{\frac{1-c}{c}}-1}\leq \frac{Cs_0^{\frac{1-c}{c}}}{2^{\frac{1-c}{c}}-1}
\end{align*}
which is independent of $y_0$. Due to the arbitrary choice of $y_0$ and the assumption $(M,g)$ is an open manifold, using Hop-Rinow theorem, we obtain  the contradiction.
Thus,  \eqref{eq:sharp-defect} holds.
\end{proof}

Combining Theorem \ref{thm:defect-bound} and Hang-Yang's identity \eqref{HY-identity}, we obtain the following corollary.
\begin{corollary}\label{cor:HY-deficit}
	Let $(M,g)$ be a complete  four-manifold with $Q_g\geq 0$ and $R_g>0$. Then
	$$	\int_MQ_g\dmu
	+\frac12\int_{M\setminus\{p\}}|\Ric_{G_p^2g}|_g^2\dmu
	=16\pi^2 \Longleftrightarrow  M\; \mathrm{is \; compact}.$$ 
\end{corollary}
\section{Proof of main theorem and further discussions}\label{sec:proof}

\begin{proof}[Proof of Theorem \ref{main theorem}]
	Since $Q_g\geq 0$ and $R_g\geq 0$, using Lemma \ref{lem:GM-lemma}, one has $\Ric_g\equiv 0$ or $R_g>0$.
	Thus,  if $g$ is not Ricci-flat, we have  $Q_g\geq 0$ and $R_g>0$.
	
First,	combining Proposition \ref{prop:finitepi} with Proposition \ref{prop:finitepiforQ=0}, we show that $\pi_1(M)$ is finite.
	
Second, with the help of Theorem \ref{thm:yamabe-total-defect}, Theorem \ref{thm:defect-bound}  and \eqref{Yamabe-range}, one has
\begin{equation}\label{CVinequalityfor-Q}
	\int_{M}Q_g\dmu\leq \frac{Y(M,[g])^2}{24}-\mathcal{B}_g\leq  \frac{Y(M,[g])^2}{48}\leq 8\pi^2.
\end{equation}
	
	Let $\pi:(\widetilde M,\tilde g)\to(M,g)$ be the universal
	cover, with $\tilde g=\pi^*g$.  Then, \eqref{CVinequalityfor-Q} also holds for $(\widetilde M,\tilde g)$. Thus, we obtain
	$$|\pi_1(M)|\int_MQ_g\ud\mu_g=\int_{\widetilde{M}}Q_{\tilde{g}}\ud\mu_{\tilde g}\leq 8\pi^2.$$
	If the equality holds, one has $Y(\widetilde M,[\tilde g])=8\sqrt 6\pi$ and $\mathcal{B}_{\tilde g}=8\pi^2.$ Consequently, using \eqref{def-B_g}, one has
	$$\int_{\widetilde M}|\Ric_{G_p^2\tilde g}|^2_{\tilde g}\ud\mu_{\tilde g}=0, \quad \forall p\in \widetilde M$$
	which yields
	$$\Ric_{G_p^2\tilde g}\equiv 0, \quad \forall p\in  \widetilde M.$$
	Choose two different points $p,q$ in $\widetilde M$.
For all $x\in \widetilde M\setminus\{p,q\}$, using Lemma \ref{lem:Brinkman} and the fact $G_p/G_q$ is not a constant, we show that $G_p^2\tilde g$ is flat. Then, using the conformal covariant property of Weyl tensor, one has
$W_{\tilde g}\equiv 0.$   Proposition \ref{prop:cylinder} yields  $(\widetilde M, \tilde g)$ is isometric to  Euclidean $(\mr^4, g_u=e^{2u}|\ud x|^2)$.

Denote $\Gamma$ as the group of  the deck
transformations which is a finite group since $\pi_1(M)$ is finite. Thus,  $\Gamma$ acts on $\mathbb R^4$ freely.
By Liouville's theorem, every element of $\Gamma$ is a Euclidean
similarity and hence an affine map. Set
$$
x_0=\frac{1}{|\Gamma|}\sum_{\gamma\in\Gamma}\gamma(0).
$$
For every $\sigma\in\Gamma$, affineness gives
$$
\sigma(x_0)
=\frac{1}{|\Gamma|}
\sum_{\gamma\in\Gamma}(\sigma\circ\gamma)(0)
=x_0.
$$
Thus $\Gamma$ has a fixed point. Since its action is free,
$\Gamma$ is trivial. Consequently, $(M,g)$ itself is isometric to $(\mr^4,  g_u)$.  Then 
\begin{equation}\label{equl-Q-integral}
	\int_{\mr^4}Q_{g_u}\ud\mu_{g_u}=\int_MQ_g\ud\mu_g=8\pi^2
\end{equation}
Applying Theorem 1.1 of  \cite{LWX} and Theorem 1.1 of  \cite{LX} together, \eqref{equl-Q-integral} yields  $I_{g_u}(\mr^4)=0$. 

Conversely, if $(\mr^4,g_u)$ has nonnegative $Q$-curvature and scalar curvature as well as $I_{g_u}(\mr^4)=0$, using  \cite[Theorem 1.1, Lemma 3.7]{LWX} and  \cite[Theorem 1.1]{LX}, we show that the integral of $Q$-curvature equal to $8\pi^2$. Since $|\pi_1(\mr^4)|=1$, we show that  the equality of \eqref{CV-type-inquality} holds.

Thus, we finish the proof.
\end{proof}

For a Riemannian  manifold $(M^n,g)$ with dimension $n\geq 3$,  Branson's  $Q$-curvature is given by 
\begin{equation}\label{eq:Q-for-n}
	Q_g
	=
	-\frac{1}{2(n-1)}\Delta_g R_g
	-\frac{2}{(n-2)^2}|\Ric_g|_g^2
	+\frac{n^3-4n^2+16n-16}
	{8(n-1)^2(n-2)^2}R_g^2
\end{equation}

In the three-dimensional case, the positivity of the $Q$-curvature and the scalar curvature gives more restrictive control over the topology of the manifolds with the help of the algebraic structure of the $Q$-curvature and the results of  \cite{CMMR, Qian}. 
For a Riemannian metric \(h\) on an \(n\)-manifold, a smooth function \(f\),
and a real number \(m>n\), define the modified  Bakry-\'Emery Ricci tensor 
\begin{equation}\label{BE-Ricci}
	\Ric^{m}_{h,f}
	=	\Ric_h+\nabla_h^2f-\frac1{m-n} \ud f\otimes  \ud f.
\end{equation}

\begin{theorem}\label{thm:three-manifolds}
	Let $(M^3,g)$ be a complete three-manifold with $Q_g\geq 0$ and $R_g\geq 0$.
	Then either $g$ is flat or $M^3$ is compact. 
\end{theorem}
\begin{proof}
	Following the same argument as in Lemma \ref{lem:GM-lemma}, one has $\Ric_g\equiv 0$ or $R_g>0$.
	In three-dimensional case, Ricci-flatness implies $g$ is flat.
	Now, we deal with the case $R_g>0.$ Let $\rho$ be the smallest eigenvalue of the Ricci tensor as in \eqref{def-of-rho}.
	Then
	\begin{equation}\label{eq:eigenvalue}
		|\Ric_g|^2\geq\rho^2+\frac12(R_g-\rho)^2.
	\end{equation}
	Using \eqref{eq:Q-for-n} and \eqref{eq:eigenvalue}, a direct  computation yields
	\begin{align*}
		-\Delta_g R_g+2\rho R_g
		&=4Q_g+8|\Ric_g|^2-\frac{23}{8}R_g^2+2\rho R_g\\
		&\geq4Q_g+12\left(\rho-\frac{R_g}{4}\right)^2+\frac{3}{8}R_g^2\geq\frac{3}{8}R_g^2.
	\end{align*}
	Since \(R_g>0\), division by \(2R_g\) gives
	\begin{equation}\label{eq:tensor}
		\Ric_g-\frac{\Delta_gR_g}{2R_g}g\geq\frac{3}{16}R_gg\geq 0.
	\end{equation}
	Then, using \eqref{eq:tensor} and  \cite[Theorem 3.13]{CMMR},  the conformal metric $\bar g=R_gg$	is complete. For brevity, set $f=\frac12\log R_g.$
	Then, one has 
	$$
	\Ric_{\bar g}=\Ric_g-\nabla^2f+\ud f\otimes \ud f
	-\left(\Delta_g f+|\nabla f|^2\right)g,$$
	and
	$$\bar\nabla^2f=\nabla^2f-2 \ud f\otimes  \ud f+|\nabla f|^2g.$$
	Thus, for the conformal metric $\bar g=R_gg=e^{2f}g$ on $M^3$ with 	$m=4$ in \eqref{BE-Ricci},  one has 
	\begin{align*}
		\Ric^{\,4}_{\bar g,f}
		&=\Ric_{\bar g}+\bar\nabla^2f- \ud f\otimes \ud f\\
		&=\Ric_g-(\Delta_g f)g-2\ud f\otimes \ud f\\
		&=\Ric_g-\frac{\Delta R_g}{2R_g}g
		+\frac12\left(
		|\nabla\log R_g|^2g-\ud\log R_g\otimes \ud\log R_g
		\right).
	\end{align*}
	Combining with  \eqref{eq:tensor}, one has 
	\begin{equation*}
		\Ric^{\,4}_{\bar g,f}\geq\frac3{16}\bar g.
	\end{equation*}
	Applying  \cite[Theorem 5]{Qian}, we show that $M^3$ is compact.
\end{proof}

	In contrast, for $n\geq 5$, even the positivity of both
the $Q$-curvature and scalar curvature does not guarantee finite generation of the fundamental group. Let $(\Sigma^2,h)$ be the two-dimensional surface with Gaussian curvature equal to $-1$.  There exist such complete hyperbolic surfaces whose fundamental groups are not finitely generated. For $n\geq 5$ and a constant $a>0$, let $(\ms^{n-2}, a^{-1}g_{\ms^{n-2}})$ be a round sphere with  sectional curvature  equal to $a$. Consider the product manifold
\begin{equation}\label{eq:product}
	(M^n,g)=(\ms^{n-2},a^{-1}g_{\ms^{n-2}})\times(\Sigma^2,h).
\end{equation}
Then, the  eigenvalues of Ricci tensor are
\[
(n-3)a,\quad \cdots,\quad (n-3)a,\quad-1,\quad-1.
\]
Using \eqref{eq:Q-for-n} and choosing  $a=20$, it is easy to check that 
\begin{equation}\label{eq:curvatures}
	R_g>0, \quad Q_g>0.
\end{equation}
Since  $\mathbb{S}^{n-2}$ is  simply connected, there holds
\[
\pi_1(M^n)\cong\pi_1(\Sigma^2).
\]
Thus, if $\pi_1(\Sigma^2)$ is not finitely generated, so is $\pi_1(M^n)$.

 If the complete manifold $(M^n,g)$ with $n\geq 5$, $Q_g\geq 0$ and $R_g\geq 0$ is conformal equivalent to a domain of a closed  manifold, the work of Ma-Qing \cite[Theorem 1.4]{MQ} shows that the    Hausdorff dimension of  the complement  is at most $(n-4)/2$. See also the foundational works \cite{CQYinv, CHY}.

Thus, the extent of the topological impact of $Q$-curvature in higher dimensions remains an open question, which we hope to explore in future work.

\vspace{3em}

\subsection*{Acknowledgement}
The author acknowledges support from the Start-up Fund
(P0064878) of The Hong Kong Polytechnic University.

\subsection*{Disclosure on AI assistance.}
During the preparation of this work, the author  used LLMs for mathematical exploration,  checking calculations and editing the exposition. The author takes full responsibility for the final manuscript.

\end{document}